\documentclass[12pt]{article}
\usepackage{fullpage}

\usepackage{graphicx,amsthm} 
\theoremstyle{definition}
\newtheorem{example}{Example}

\usepackage[dvipsnames]{xcolor}

\usepackage{caption}
\usepackage{subcaption}
\usepackage{comment}
\usepackage{authblk}
\usepackage[numbers,sort&compress]{natbib}
\usepackage{doi}

\usepackage{appendix}

\usepackage{graphicx}
\graphicspath{ {./FiguresSST/} }
\usepackage{float} 
\usepackage[export]{adjustbox} 

\usepackage{breqn}
\usepackage{titlesec}

\usepackage{bm}
\usepackage{dsfont}

\usepackage{amsmath,amsfonts,amsthm,bm} 
\usepackage{upgreek}
\usepackage{mathtools}

\usepackage{amssymb}

\usepackage{amsfonts}

\usepackage{geometry}
\usepackage{tabularray}

\usepackage{sectsty}

\usepackage{tikz}
\usetikzlibrary{shapes.geometric, arrows}
\usetikzlibrary{decorations.pathmorphing,shapes}
\usepackage{makecell}
\usepackage{booktabs,caption}
\usepackage[flushleft]{threeparttable}

\usepackage{nomencl}
\makenomenclature

\usepackage{color}   
\usepackage{hyperref}
\hypersetup{
    colorlinks=true, 
    linktoc=all,     
    linkcolor=blue,  
}

\usepackage{framed}

\usepackage{matlab-prettifier}

\title{Data-driven methods for computation of \\ optimal linear response in \\ high-dimensional dynamical systems}

\author[1]{Gary Froyland}
\author[2]{Dimitrios Giannakis}
\author[1]{Nicholas Peters}

\affil[1]{School of Mathematics and Statistics, UNSW Sydney\authorcr Sydney NSW 2052, Australia}
\affil[2]{Department of Mathematics, Dartmouth College\authorcr Hanover, NH 03755, U.S.A.}

\begin{document} 

\maketitle

\begin{abstract}
    We develop a data-driven framework for estimating optimal linear response of nonlinear dynamical systems.
    Our approach is based on kernel-smoothed approximations of the transfer/Koopman operators of the system, built from possibly high-dimensional observations along trajectories.
    Combining these operator approximations with the theory developed in [Antown et al.\ (2018), \emph{J.\ Stat.\ Phys.}, 170(6), 1051--1087], we formulate a computationally tractable optimization problem for the optimal infinitesimal perturbation realising any desired manipulation of the dynamics.
    We also introduce a notion of optimal-response vector fields for visualising, and physically interpreting, the system's response to the optimal perturbation under arbitrary observations.
    
    While our data-driven techniques apply to general optimal linear response problems, the focus in this current work is on finding perturbations of the spectrum that optimally increase the frequency or optimally suppress the decay of correlations of almost-cycles or almost-invariant sets associated with the eigenvalues of the kernel-smoothed transfer operator.
    We illustrate our approach with applications to low-dimensional periodic and chaotic systems, as well as a high-dimensional example involving the El Ni\~no Southern Oscillation in a comprehensive Earth system model.
    In these examples our approach discovers nontrivial optimal perturbations of the system, which are \emph{post hoc} natural and consistent with the desired dynamical objectives.
\end{abstract}

\section{Introduction}

Estimating the response of dynamical systems to external perturbations is an important problem across the spectrum of science and engineering disciplines dealing with complex systems, including fluids, granular matter, and biomolecules.
Frequently, the system under study is not accessible to direct perturbative study, and one is faced with the problem of estimating the response from observations of the unperturbed system.
A prototypical example is the Earth's climate system where we have observational access to only a single realization of the dynamics.
There is considerable interest in estimating changes in the statistical behavior of the system under natural and anthropogenic perturbations such as volcanoes and greenhouse gas emissions, respectively \cite{GhilLucarini20, LucariniEtAl17}.     
Linear response theory is a framework for addressing this problem when the perturbation is sufficiently weak and the statistical response depends smoothly (differentiably) on it in a suitable sense \cite{Ruelle97, Ruelle09, Baladi14}.

The modern formulation of the theory has its roots in various types of fluctuation--dissipation theorems (FDTs), which link the linear response of observables to certain correlation functions computed for a system in thermodynamic equilibrium. 
Early versions of the FDT were put forward by Callen, Welton, Kubo and other authors in the 1950s \cite{CallenWelton51,Kubo66} for linear systems under Langevin dynamics, generalizing previously known relations for Brownian motion and Johnson noise.   
Kraichnan \cite{Kraichnan59} later extended this work to general (nonlinear) Hamiltonian systems, again in thermodynamic equilibrium. 
In work from the 1970s \cite{Leith75}, Leith studied versions of the FDT in geophysical fluid flow, and gave a justification for the applicability of the theorem in certain problems for observables of large-scale atmospheric dynamics, despite the fact that these systems are neither Hamiltonian nor are in thermodynamic equilibrium.   

Later on, in the mid 1990s, Gallavotti and Cohen \cite{GallavottiCohen95,Gallavotti96} developed a formulation of the FDT for nonequilibrium systems using the framework of modern dynamical systems theory.  
The class of systems they considered is that of dissipative, hyperbolic, chaotic systems possessing physical or Sinai-Ruelle-Bowen (SRB) invariant probability measures \cite{Young02}.
These invariant measures do not, in general, require thermodynamic equilibrium, yet they exhibit several useful properties for linear response, including differentiability of response functionals and approximability of the statistics of the invariant measure from time averages of observables along dynamical trajectories.
They demonstrated validity of the FDT in such settings using numerical experiments.

In more detail, classical linear response \cite{Ruelle97,Baladi14} for discrete-time dynamical systems $T:M\to M$ is frequently concerned with the derivative (with respect to a system parameter) of a physical probability measure for a map $T$ with some hyperbolicity properties on a manifold $M$.
Hyperbolicity in the state space of the dynamics (the manifold $M$), means that at each $x\in M$, the local linearised action of the dynamics of $T$ acting on the tangent space $\mathcal{T}_xM$ of $M$ at $x$ is a combination of strict expansion and strict contraction, occurring in complementary subspaces of $\mathcal{T}_xM$.
For such dynamics, under further mild assumptions, it is known \cite{BowenRuelle75} that there is a unique probability measure on $M$, called the SRB measure, which is describes the long-term distribution of trajectories initialised at  Lebesgue almost-every $x\in M$.
A natural question is ``how does this SRB measure change under small perturbations to the dynamics?''.
A rigorous proof of linear response  in this uniformly hyperbolic class was given by Ruelle \cite{Ruelle97}. Since then a vast amount of linear response theory has been developed for different classes of deterministic and stochastic dynamical systems including piecewise-expanding, intermittent, nonuniformly hyperbolic, stochastic, and high-dimensional chaotic systems \cite{GouezelLiverani06,AbramovMajda08,Ruelle09,BaladiSmania09,HairerMajda10,Baladi14,BaladiKunaLucarini17, BaladiTodd16, WormellGottwald18}.
A wide variety of numerical approaches to computing response have also been proposed, e.g.~\cite{GottwaldWormellWouters16,bahsoun2018rigorous,  ChandramoorthyWang22,ni2026fast,FroylandPhalempin25}.
In climate dynamics, a large body of literature \cite{Bell80,Gritsoun01,GritsounEtAl02,MajdaEtAl05,MajdaEtAl10b,LucariniEtAl17} has explored numerical methods and applications of linear response theory building on Leith's work.  

In parallel with the development of this body of work on linear response, starting from the late 1990s, there has been considerable interest in the development of data-driven techniques for approximation of Koopman and transfer operators of dynamical systems; see e.g.\ \cite{DellnitzJunge99,SchutteEtAl01,Mezic05,MezicBanaszuk04} for some of the earlier references and \cite{OttoRowley21,Colbrook24,BruntonEtAl22} for reviews.
From an analytical standpoint, formulations of linear response have strong connections with operator-theoretic constructions.
For example, operator-theoretic formulations of response theory express derivatives of statistical quantities through perturbations of transfer operators acting on suitable anisotropic Banach spaces, or dually through Koopman operators on observables \cite{Ruelle97, GouezelLiverani06, Baladi2018}.

Recent work \cite{SantosGutierrezLucarini22,ZagliEtAl26,LucariniEtAl26} has leveraged these connections to build computational methods for linear response in stochastic systems that take advantage of mature data-driven Koopman and transfer operator techniques such as the extended dynamic mode decomposition (EDMD) \cite{WilliamsEtAl15}. In this approach, eigendecompositions of Koopman operators are used to estimate correlation functions and other statistical quantities appearing in response formulas, bypassing some of the challenges associated with direct computation methods.

The focus of the present work is on the linear response of the outer spectrum of the transfer operator, which encodes the dynamics of slowly decaying features such as almost-invariant sets \cite{DellnitzJunge99, DellnitzEtAl00} and almost-cyclic sets \cite{DellnitzJunge99,FroylandEtAl14b}.  
In the uniformly hyperbolic setting,
one considers a family of maps $\{T_\delta\}_{\delta\in [0,\bar\delta]}$, and writes $T_\delta=T_0+\delta\cdot \dot{T}+o(\delta)$, where $\dot T \colon M\to\mathbb{R}^d$ is a smooth vector field. Denoting by $\mathcal{L}_\delta$ the transfer operator\footnote{considered as acting on a suitable Banach space \cite{GouezelLiverani06}.} of $T_\delta$ with respect to Lebesgue measure on $M$, one can rigorously derive \cite{GouezelLiverani06,  Porte19, FroylandPhalempin25} the first-order Taylor expansion $\mathcal{L}_\delta=\mathcal{L}_0+\delta\cdot \dot{\mathcal{L}}+o(\delta)$, where the operator $\dot{\mathcal{L}}$ represents the derivative of $\mathcal{L}_\delta$ with respect to $\delta$. Each $\mathcal{L}_\delta$ has its own outer spectrum whose derivatives are guaranteed by e.g.\ \cite{GouezelLiverani06}. 

\subsection{Our contributions}

We study the problem of \emph{optimal} linear response for the outer spectrum of the transfer operator; that is, finding a linear perturbation of a dynamical system that maximally affects the dynamics of slowly decaying features. 
The strength of invariance of almost-invariant sets is quantified by the magnitude of the corresponding real eigenvalue of the transfer operator \cite{DellnitzJunge99,DellnitzEtAl00,Froyland05}, and similarly the decorrelation rate of almost-cycles are quantified by the magnitude of the corresponding complex eigenvalues \cite{DellnitzJunge99,FroylandEtAl14b}.
Furthermore, the period of almost-cycles are given by the argument of the corresponding eigenvalue \cite{DellnitzJunge99,FroylandEtAl21,CastroFroyland25}.
By optimally manipulating the magnitude and argument of real or complex eigenvalues one may control how well almost-invariant sets trap trajectories, how rapidly almost-cycles decorrelate, and the period of almost-cycles.
The optimal $\dot{T}$ tells one how to make an infinitesimal perturbation to the original dynamics $T$ to achieve the greatest such effects post-perturbation.

In the context of enhancing fluid mixing, optimal linear response approaches have been used to maximally manipulate the second-largest eigenvalue of the evolution (Perron--Frobenius) operator for the flow in periodically forced \cite{froyland2017optimal} and aperiodically forced \cite{froyland2020computation} flows.
Closer to the present work, for finite-state Markov chains Antown \emph{et al.}~\cite{AntownFroyland18} determined perturbations optimising the location of the second largest eigenvalue (and also perturbations to maximise the change in the equilibrium distribution and maximise the expectation of an observation function, both for non-stationary sequences of Markov chains).
The optimal spectral manipulation results for finite-state Markov chains were extended to the large and flexible class of Markov Hilbert--Schmidt operators in \cite{AntownEtAl22}. 
In recent work, Santos Guti{\'e}rez \emph{et al.}\ \cite{SantosGutierezEtAl25} applied the methods of \cite{AntownFroyland18} to identify optimal perturbations of various entropic functionals, using Markov chain estimates of the transfer operator derived from Ulam's method and a collection of unstable periodic orbits.

The theory and numerical methods of \cite{AntownFroyland18,AntownEtAl22} that determine optimal spectrum-manipulating perturbations, require the dynamical flow or map to be known, or at least to be able to access some numerical oracle on the flow or map. 
In this work our goal is to extend the framework of \cite{AntownFroyland18,AntownEtAl22} in two key ways. 
Firstly to the situation where the dynamical system is unknown, and one only has a (possibly single, as in the case of historical climate data) time series of typically non-injective observations. 
Secondly to the setting where the observations lie in some high-dimensional space, precluding the use of explicit meshing or partitioning approaches such as Ulam's method.
To address these computational issues we use kernel smoothing methods for operator approximation \cite{DasGiannakis19, Giannakis21a}, optionally in conjunction with delay-coordinate maps \cite{Takens81}.
The use of kernels restores Markovianity of the ``deterministic'' observed data under partial observations, and enables the construction of data-driven approximations of the Koopman operator.

We then adjust the framework of \cite{AntownFroyland18} to compute optimal perturbations of the Markov chain, which maximally alter the magnitude or argument of a selected real or complex eigenvalue of the matrix representation of the Markov chain (and of the Koopman operator).
In a final crucial contribution for applications, we devise a computational procedure to reconstruct vector fields associated with the optimal Markov chain perturbation in arbitrary data spaces. In this procedure, an optimal-response vector field is reconstructed for any target observable of interest, thus providing a physically interpretable visualization of the optimal perturbation.
While our focus in this work is on optimising almost-cycles and almost-invariant sets through the linear response of the transfer/Koopman operator spectrum, our data-driven framework may readily be applied to general optimal linear response problems, including the classical linear response discussed in the introduction.

We present applications of our approach in a suite of examples of increasing degrees of complexity: (i) uniform and non-uniform circle rotations; (ii) the Lorenz 63 chaotic system; and (iii) the El Ni\~no Southern Oscillation (ENSO) extracted from simulated Indo-Pacific sea surface temperature (SST) evolution in a comprehensive Earth System model---the Community Climate System Model Version 4 (CCSM4). We compute perturbations that either (i) optimally suppress decay of correlations of these (almost)-cycles or (ii) increase the speed/frequency of the cycles.

Some sample findings are as follows. Our analysis of the non-uniform circle rotation shows that in order to maximally speed up the overall cycle, one should perturb tangentially to the current rotation direction, but in ``inverse proportion'' to the current speed.
For the Lorenz flow, studying the symmetric main rotation around the two wings of the Lorenz attractor, we find that a related strategy arises in the perturbation that maximises an increase in wing-rotation speed. One perturbs orthogonal to level sets of the complex eigenvector with strength inversely proportional to the original speed in the same direction. Our analysis of the two dominant almost-invariant sets in the Lorenz attractor reveals a relatively simple optimal strategy to maximally perturb to further strengthen the almost-invariance.
{In the ENSO example we use a time series of Indo-Pacific SST fields to perform our optimisation.
We then reconstruct the optimal perturbation vector field in terms of both SST and surface wind fields.
Firstly, we find that the optimal ENSO frequency-increasing perturbation leads the ENSO cycle, as identified by the Koopman/transfer operator eigenfunctions from \cite{FroylandEtAl21}, by a phase of $\pi/2$.
Secondly, the perturbation that optimally suppresses decay of correlations is largely in phase with the cycle.}
{Phase composites of these optimal-response vector fields, conditioned on the phase of the ENSO eigenfunction, exhibit significant activity in the tropical Indian Ocean.
Our results suggest a specific connection between the Indian Ocean dipole (IOD) \cite{HameedEtAl99,WebsterEtAl99} phasing and these predictability-enhancing perturbations.
}

The paper is organized as follows.
In section~\ref{sec:operator_approximation} we describe our kernel-based approach for Koopman operator approximation and illustrate the interpretation of its spectrum for a circle rotation.
In section~\ref{sec:optlinresp} we describe how we infer dynamical properties from the Koopman spectrum, setup the various optimisation problems, and describe a novel approach to infer perturbations in the phase space of the dynamics from perturbations of the Markov chain approximation of the Koopman operator.
All constructions are illustrated for a circle rotation.
Sections~\ref{sec:fast_slow_circle}, \ref{section:lorenz}, and \ref{sec:enso} present our case studies of these techniques applied to a non-uniform circle rotation, the Lorenz system, and ENSO, respectively.
We conclude in section~\ref{sec:conclusions}.

\section{Koopman operator approximation}
\label{sec:operator_approximation}

\subsection{Observable functions and time delay-embedding}

Let the manifold $M$ be the state space of a dynamical system and $\Phi^t: M \rightarrow M, t \in \mathbb{R}$, a continuous flow on $M$ with invariant probability measure $\mu$. 
Let $\omega_0$ be some initial point and define $\omega_i = \Phi^{i \Delta t}(\omega_0)$, so that $\{\omega_i\}_{i=0}^N \subset M$, is the discrete-time state space trajectory with sampling interval $\Delta t \in \mathbb{R}$ .
The state space $M$ may be unknown, in the sense that we cannot obtain complete knowledge of any given state $\omega \in M$. A common situation is where it is possible to observe $x = F(\omega) \in \mathbb{R}^d$ with an observation function $F: M \rightarrow \mathbb{R}^d$. Measuring with $F$ at a fixed sampling interval $\Delta t$ for $N+1$ time steps produces the trajectory of observations $\{x_i\}_{i=0}^N \in \mathbb{R}^d$ where $x_i = F(\omega_i)$.

Often the observation $x_i$ does not fully capture all information about the state $\omega_i$.
That is, the observation map $F$ is not injective and there is no dynamical evolution map $\Psi^t \colon \mathbb R^d \to \mathbb R^d$ giving the evolution of the observations as $x_i = \Psi^{i\,\Delta t}(x_0)$.
For example, in the ENSO case study in section~\ref{sec:enso} below $F$ will be the map that sends the entire state $\omega$ of the CCSM4 model (consisting of three-dimensional oceanic and atmospheric fields, as well as a multitude of other variables representing different Earth system processes) to SST fields sampled at $d$ gridpoints in an Indo-Pacific domain.
In this case, a dynamical map $\Psi^t$ on observations space does not exist since the evolution of Indo-Pacific SST is influenced by other state variables such as atmospheric winds and sub-surface oceanic circulation.

To improve our representations of the states in such partial observations scenarios, a classical approach, which we employ here, is the Takens delay-embedding technique \cite{Takens81,SauerEtAl91}.
The idea behind this method is that even though individual snapshots, $x_i$, may not contain sufficient information to determine the future evolution, \emph{sequences} of samples $x_i, x_{i+1}, \ldots$, implicitly contain sufficient information from $\Phi^t$ so as to uniquely determine their future evolution.
Returning to our Indo-Pacific SST example, it is known from geophysical fluid dynamics that the time history of the evolution of the surface ocean carries imprints of subsurface oceanic circulation or atmospheric wind forcing, which influences the future evolution of SST.
This makes the evolution of SST sequences more predictable than individual snapshots.

For a chosen parameter pair $Q,L \in \mathbb N$, representing the number of delays and delay length, respectively, we define the delay embedding map $F_{Q,L}: M \rightarrow \mathbb{R}^{\tilde{d}}$ as $F_{Q,L}(\omega) = [F(\omega), F(\Phi^{L \Delta t}\omega), \ldots, F(\Phi^{(Q-1)L \Delta t} \omega)]$, where $\tilde{d} = Qd$ is the embedding dimension. Note that the values of $F_{Q,L}$ on the states $\omega_i$ are given by concatenation of $Q$ observation vectors $x_{i+qL} \in \mathbb{R}^d, q = 0,\dots, Q-1$, each separated by a delay of $L$, namely $F_{Q,L}(\omega_i) = [x_i, x_{i+L}, x_{i+2L}, \dots, x_{i+(Q-1)L}]$. In particular, the values $\tilde x_i := F_{Q,L}(\omega_i)$ can be obtained without explicit knowledge of $\omega_i$. The delay-embedded trajectory is $\{\tilde{x}_i\}_{i=0}^{\tilde{N}}$ with $\tilde{x}_i = F_{Q,L}(x_i)$, where $\tilde{N} = N -(Q-1)L$ is the length after the truncation from delay embedding. When delay-embedding is unnecessary, we set $Q = 1$ so that $\tilde{x}_i = x_i$ for any choice of $L$. Various delay-embedding theorems established over the past four decades (e.g., \cite{Takens81,SauerEtAl91,Robinson05,Robinson08,DeyleSugihara11}) show that under appropriate mild assumptions on the dynamical system $\Phi^t$, observation map $F$, and sampling interval $\Delta t$, for sufficiently large $Q$ and fixed $t$ there exists a map $\tilde{\Phi}^t:\mathbb{R}^{\tilde{d}} \rightarrow \mathbb{R}^{\tilde{d}}$ that is compatible with the state space dynamics on $M$, in the sense of the conjugacy relation $\tilde \Phi^t \circ F_{Q,L} = F_{Q,L} \circ \Phi^t$. This implies that $\tilde \Phi^t$ takes each $\tilde{x}_i$ to its successor $\tilde{\Phi}(\tilde{x}_i) = \tilde{x}_{i+1}$. Likewise, the embedding space $\mathbb{R}^{\tilde d}$ inherits an invariant measure $\tilde{\mu}=\mu\circ (F^{Q,L})^{-1}$ from the invariant measure $\mu$ on the state space. 

\subsection{A kernel-based Markov transition matrix approximation for the transfer/Koopman operator}

We will use the observed trajectory to approximate the transfer operator $\mathcal{L}^t f = f \circ \tilde{\Phi}^{-t}$ and its adjoint, the Koopman operator $\mathcal{K}^t f = f \circ \tilde{\Phi}^t$, defined on functions $f: \mathbb{R}^{\tilde{d}} \rightarrow \mathbb{R}^{\tilde{d}}$ in $L^2(\tilde\mu)$. 
We build kernel-smoothed approximations of these operators using the trajectory data $\{\tilde x_i \}_{i=0}^{\tilde N}$.
To this end, for $\epsilon > 0$, define the Gaussian kernel $k_\epsilon:\mathbb{R}^{\tilde{d}}\times \mathbb{R}^{\tilde{d}}\to \mathbb{R}^+$ by
\begin{equation}
\label{eq:kernel}    k_\epsilon(\tilde{x},\tilde{y}) = \exp\left(-\frac{\|\tilde{x}-\tilde{y}\|^2}{\epsilon^2 \tilde{D}^2}\right),
\end{equation}
where $\tilde{D} = \frac{1}{(\tilde{N}+1)^2} \sum_{i,j = 0}^{\tilde{N}} \|\tilde{x}_i-\tilde{x}_j\|$ is the mean pairwise distance between the points $\{\tilde{x}_i\}_{i=0}^{\tilde{N}}$ in the embedding space\footnote{We use the $\tilde{D}$ term in~\eqref{eq:kernel} to make the scaled distances $\frac{\|\tilde{x}-\tilde{y}\|}{\tilde{D}}$ less dependent on $Q$ and $L$, so that for a given observation time series, the choice of $\epsilon$ remains relatively stable.}. Following \cite{FroylandEtAl21}, we define the operator $\mathcal K^t_\epsilon$ on $L^2(\mu)$ as
\begin{equation}
    \label{eq:markov_op}
    \mathcal{K}^t_\epsilon f(\tilde{x}) = \int_{\mathbb{R}^{\tilde{d}}} \frac{k_\epsilon(\tilde{x},\tilde\Phi^{-t}(\tilde{y}))}{\int_{\mathbb{R}^{\tilde{d}}} k_\epsilon(\tilde{x}, \tilde\Phi^{-t}(\tilde{z})) \, d\tilde{\mu}(\tilde z)} f(\tilde y) \, d\tilde{\mu}(\tilde y).
\end{equation}
Note that the expression \eqref{eq:markov_op} is the same as the expression for the approximation $P_\epsilon f$ for the transfer operator given in \cite{FroylandEtAl21}, except that in \eqref{eq:markov_op} the $\Phi$ has been replaced with $\Phi^{-1}$ because here we approximate the Koopman operator instead of the transfer operator.

One readily verifies that $\mathcal K^t_\epsilon$ maps positive functions to positive functions and it preserves constant functions.
{Rewriting \eqref{eq:markov_op} as 
\begin{equation}
    \label{eq:markov_op2}
    \mathcal{K}^t_\epsilon f(\tilde{x}) = \int_{\mathbb{R}^{\tilde{d}}} p_\epsilon(\tilde x, \tilde y) f(\Phi^{t}(\tilde y)) \, d\tilde{\mu}(\tilde y), \quad p_\epsilon(\tilde x, \tilde y) := \frac{k_\epsilon(\tilde{x},\tilde{y})}{\int_{\mathbb{R}^{\tilde{d}}} k_\epsilon(\tilde{x}, \tilde{z}) \, d\tilde{\mu}(\tilde z)},
\end{equation}
we see that on the right-hand side we have a smoothed approximation -- reminiscent of the Naradaya--Watson kernel regression \cite{Nadaraya64,Watson64} -- of $f \circ \tilde \Phi^t$ based on the normalized kernel $p_\epsilon$. In the $\epsilon\to 0$ limit, general sufficient conditions for $\mathcal K^t_\epsilon f(\tilde x)$ to converge to $\mathcal Kf(\tilde x)$ are that: integration by $p_\epsilon(\tilde x, \cdot)$ is an approximate identity (to localise the averaging) and that $\tilde x$ is a $\tilde\mu$-Lebesgue point of the integrand (to ensure that a local average is meaningful).
This result can be extended to convergence in $L^p$ norm and almost everywhere convergence using the classical theory of Hardy--Littlewood maximal functions \cite{Stein93}, coupled with conditions on measure density (or ``thickness'') that prevent the normalization function $\tilde x \mapsto \int_{\mathbb R^d} k_\epsilon(\tilde x, \tilde z) \, d\tilde\mu(\tilde z)$ from vanishing as $\epsilon \to 0$; e.g., \cite{hajlasz_sobolev_2008,GiannakisLatifiJebelli26}. Variants of kernel smoothing theory have also been developed when $\tilde\mu$ is not absolutely continuous with respect to the Lebesgue measure; e.g., under appropriate capacity conditions on metric measure spaces \cite{BjornEtAl01,CantoEtAl25}. In the context of dynamical systems (as opposed to general function approximation problems) operators of the form \eqref{eq:markov_op2} have been studied as ``small random perturbations'' of dynamical systems \cite{Khasminskii63, Kifer74}, particularly by Kifer, and have found application in numerical approximations of transfer operators by Ulam's method \cite{Ulam64}, where the kernels are indicators on fine partitions of the domain of $\Phi^t$.} If $\epsilon>0$ then rather than sending $f(\tilde{x})$ directly to $f(\tilde{\Phi}^{t}(\tilde{x}))$, expressions \eqref{eq:markov_op} and \eqref{eq:markov_op2} blur by integrating $f$ over $\mathbb{R}^{\tilde{d}}$, weighted by the distance of each point to $\tilde{\Phi}^{t}(\tilde{x})$. 

We now describe how to construct an estimator of $\mathcal{K}_\epsilon^{\Delta t}$, i.e.\ the operator $\mathcal K_\epsilon^t$ from~\eqref{eq:markov_op} for our sampling interval $t = \Delta t$.
We assume that the sampling measure $\tilde{\mu}_N=\frac{1}{\tilde N}\sum_{i=1}^{\tilde N} \delta_{\tilde x_i}$ associated with our data $\tilde x_i$ converges weakly to the invariant measure $\tilde \mu$, meaning $\frac{1}{\tilde N} \sum_{i=1}^{\tilde N} f(\tilde x_i)$ converges to $\int_{\mathbb R^{\tilde d}} f(\tilde x) \, d\tilde\mu(\tilde x)$ as $\tilde N\to \infty$ for $f$ continuous.
Numerically, we know only the $\tilde{N}+1$ points $\{\tilde{x}_i\}_{i=0}^{\tilde{N}}$ with $\tilde{x}_{i+1} = \tilde{\Phi}^{\Delta t}(\tilde{x}_i)$, so we construct a discrete version of $\mathcal K^{\Delta t}_\epsilon$ by replacing the integrals over $\mathbb{R}^{\tilde{d}}$ in \eqref{eq:markov_op} by time averages over the trajectory $\{\tilde{x}_i\}_{i=0}^{\tilde{N}-1} \in \mathbb{R}^{\tilde{d}}$. The result is the $\tilde{N}\times\tilde{N}$ row-stochastic matrix $P_0$, which acts on $f(\tilde{x}_i)$ via
\begin{equation*}
    P_0 f(\tilde x_i) 
    = \sum_{j=1}^{\tilde{N}} P_{0,ij} f(\tilde{x}_j),
\end{equation*}
where
\begin{equation}
    \label{eq:markov_mat}
    P_{0,ij} =  \frac{k_\epsilon (\tilde{x}_{i},\tilde{x}_{j-1})}{\sum_{j'=1}^{\tilde{N}} k_\epsilon(\tilde{x}_{i}, \tilde{x}_{j'-1})} \in [0,1]
\end{equation}
is the $(i,j)$th entry of $P_0$ for $i,j = 1,\dots,\tilde{N}$.

{Note that expression on the RHS of \eqref{eq:markov_mat} is the classical Naradaya--Watson estimate of the conditional expectation of the time-evolved function $f \circ \tilde \Phi^{\Delta t}$ given the observations under $F_{Q,L}$. More specifically, by using} \eqref{eq:markov_op}, and the fact that $\tilde{\mu}_N\to \tilde{\mu}$ weakly, we have $ P_0 f(\tilde x_i) \approx \mathcal K^{\Delta t}_\epsilon f(\tilde x_i) \approx \mathcal K^{\Delta t} f(\tilde x_i)$.
Using the fact that $k_\epsilon$ is symmetric, it is not hard to show that $P_0$ in \eqref{eq:markov_mat} is the transpose of the column-stochastic matrix used to approximate the transfer operator in \cite{FroylandEtAl21}.
Both here and in \cite{FroylandEtAl21}, we multiply vectors on the right of the matrix $P_0$;  the fact that here we use the transpose of the $P$ in \cite{FroylandEtAl21} is consistent with the fact that \eqref{eq:markov_mat} approximates the Koopman operator rather than the transfer operator. By classical results \cite{kemenysnell}, the spectrum of $P_0$ is confined to the unit disk in the complex plane.

\begin{example} \label{sect:circlerotation}

\textbf{Operator analysis of a circle rotation.} Consider a circle rotation, $\Phi^{\Delta t}\colon S^1 \to S^1$, $\Phi^{\Delta t}(\theta) = \theta + \alpha \Delta t \mod 2 \pi$, where $\alpha = 4$ is a rotation frequency parameter and $\Delta t = 0.1$ is the sampling interval. With these values, the discrete-time map $\Phi^{\Delta t}$ describes an irrational rotation by the angle $\Delta \theta = \alpha\, \Delta t = 0.4$ radians, corresponding to a period $2 \pi / \alpha = \pi/2  \approx 1.57$, or a rotation frequency $\alpha / 2\pi \approx 0.6366$.
The system is observed under the map $F\colon S^1 \to \mathbb R^2$, defined to be the standard embedding of the circle into $\mathbb R^2$, $F(\theta) = [\cos\theta, \sin \theta]$. 
We collect $N+1 = 1000$ samples $x_i = F(\theta_i) \in \mathbb R^2$ with $\theta_i = \Phi^{i\,\Delta t}(\theta_0)$ and $i = 0, \ldots, N$, on a trajectory starting from an arbitrary initial condition $\theta_0 \in S^1$.
The $N+1$ observations correspond to approximately 63.7 laps of the circle. Since time delay-embedding is unnecessary for this system, we set $Q=1$ so that $\tilde{x}_i = x_i$ for all $i = 0, \dots, \tilde{N}$, where $\tilde{N} = N$.

Using this data, we build the Markov transition matrix $P_0 \in \mathbb R^{N \times N}$ given by~\eqref{eq:markov_mat} for the kernel bandwidth parameter $\epsilon = 0.32$.
Note that this is a relatively large $\epsilon$ for such a simple trajectory, but this choice serves to better illustrate the mechanisms of the transition matrix $P_0$, shown in Figure \ref{fig:circle_rotation_unpert}(left) by increasing the width of the diagonal bands. In Figure \ref{fig:circle_rotation_unpert}(left), these diagonal bands are caused by the circle trajectory returning nearby itself on each lap of the circle.

\begin{figure}
    \centering
    \includegraphics[height=5.5cm]{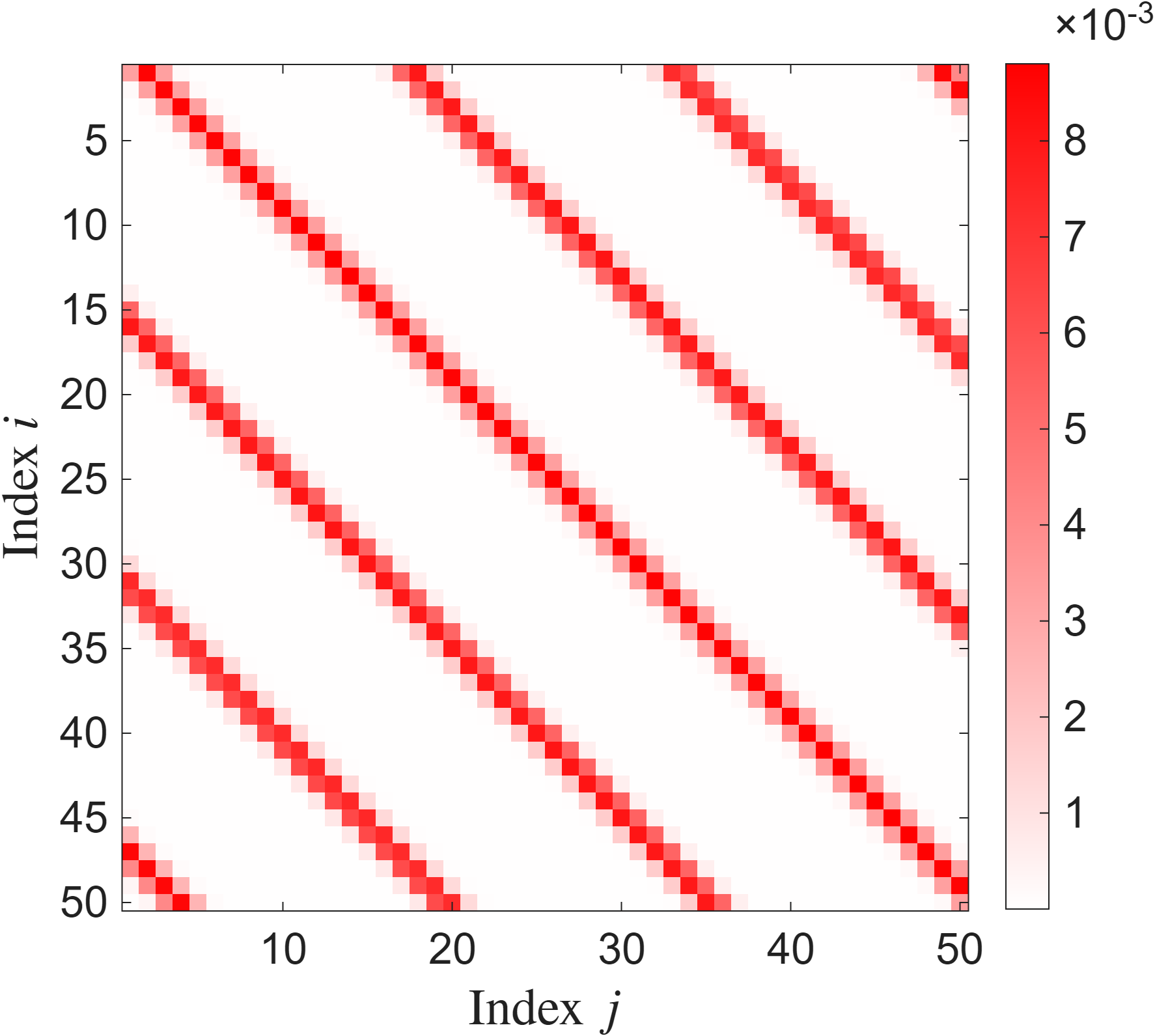}
   \includegraphics[height=5.5cm]{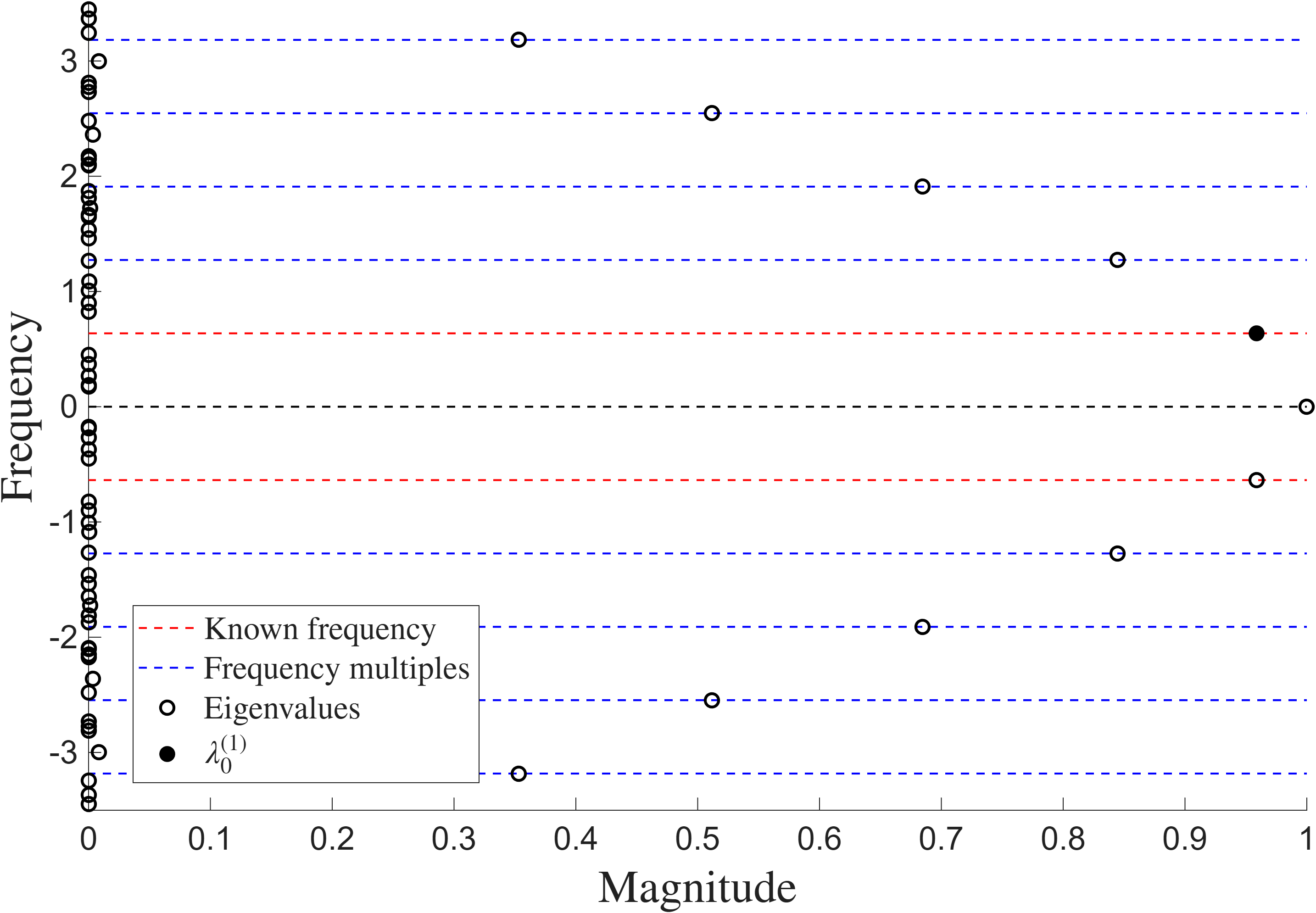}
\caption{\emph{Left:} Heat map of a $50 \times 50$ diagonal block of the $1000 \times 1000$ transition matrix $P_0$. 
\emph{Right:} Spectrum of $P_0$, visualized in terms of the magnitudes, $\lvert\lambda^{(k)}\rvert$, and frequencies $\nu^{(\pm k)} = \arg \lambda^{(\pm k)} / 2\pi \Delta t$, associated with the eigenvalues $\lambda^{(\pm k)}$.}
\label{fig:circle_rotation_unpert}
\end{figure}

Next, we compute an eigendecomposition $P_0 v^{(k)} = \lambda^{(k)} v^{(k)}$ of $P_0$, ordering the eigenvalues  in order of decreasing modulus, $1 = \lambda^{(0)} > |\lambda^{(\pm 1)}| > \cdots$, and using $(\cdot)^{(\pm\cdot)}$ to indicate complex-conjugate eigenvalue/eigenvector pairs.
The leading nontrivial eigenvalue $\lambda^{(+1)}$ is associated with the rotation of the circle and has $\arg{\lambda^{(+1)}} \approx \alpha \, \Delta t  = 0.4$, with an approximation error of $10^{-6}$. Equivalently, the frequency of the cycle associated with this eigenvalue can be written $\nu^{(+1)} = \arg{\lambda^{(+1)}}/2\pi \, \Delta t \approx \alpha/2\pi \approx 0.6366$, which agrees with the frequency of the circle rotation. This can be seen in Figure \ref{fig:circle_rotation_unpert}(right), where the second eigenvalue lies on the horizontal line representing $\alpha/2\pi$.
The leading few moduli, $|\lambda^{(\pm k)}|$, and frequencies, $\nu^{(\pm k)} = \arg \lambda^{(\pm k)} / 2\pi\, \Delta t$, are displayed in Figure \ref{fig:circle_rotation_unpert}(right).
 Numerically, the results in Figure \ref{fig:circle_rotation_unpert}(right) are in good agreement with the relationships
\begin{equation}
\label{evalbehaviour}
    \nu^{(\pm k)} = \pm k \nu^{(\pm 1)} = \pm k\alpha/2\pi, \quad |\lambda^{(\pm k)}| = 1-k^2\epsilon^2\Delta t,
\end{equation}
which may be derived for the spectrum of the Fokker--Planck operator $L_\epsilon f = - \alpha f' + \epsilon^2 f''$
for a stochastic perturbation of the circle rotation as follows.
The $k$th-order eigenfunctions are $f_k(x)=\exp(ikx)$, which have eigenvalues $-\alpha k i - \epsilon^2 k^2$.
Evolving for a time $\Delta t$, the eigenvalues of the evolution operator  $\exp(\Delta t\cdot L_\epsilon)$ become $\exp(\Delta t(-\alpha k i -\epsilon^2k^2))\approx (1-\Delta t\cdot\epsilon^2k^2)\exp(-\Delta t \alpha ki)$, for small $\Delta t, \epsilon$ and modest $k$.
Thus, we see the magnitude decays with $k$ approximately as $1-\Delta t\cdot\epsilon^2k^2$ and the argument grows as $k\alpha\Delta t$.
This argument growth corresponds to a frequency growth of $\nu_k = k \alpha \Delta t /2 \pi \Delta t = k \alpha/2\pi$ 
Thus we obtain the relationships in \eqref{evalbehaviour}.
\end{example}

\section{Optimal linear response}
\label{sec:optlinresp}

We consider perturbations of the form $P_\delta=P_0+\delta\cdot\dot P$, where the rows of the derivative matrix $\dot P$ sum to zero to ensure that $P_\delta$ is a stochastic transition matrix.
We emphasise that the first-order-in-$\delta$ expansion of $P_\delta$ involves $\dot{P}$, which is a derivative acting on the tangent space of the manifold of stochastic matrices anchored at $P_0$.
The derivative $\dot{P}$ is the discretised analogue of the adjoint of the Frech\'et derivative of the transfer operator $\dot{\mathcal{L}}$ mentioned in the introduction, and is not the same as the state-space derivative $\dot{\Phi^1}$ in the first-order expansion $\Phi^1_\delta = \Phi^1_0 + \delta\cdot \dot{\Phi^1}$.

\subsection{Eigenvalue information}

Previous work on almost-invariant and almost-cyclic sets \cite{DellnitzJunge99,  DellnitzEtAl00, Froyland05, FroylandEtAl14b} has demonstrated that  the eigenvalues close to the unit circle carry important information on the dynamics of $T_0$. Real positive eigenvalues close to 1 indicate the existence of phase space decompositions that are close to invariant under the true dynamics $\Phi^1_\delta$.
Such decompositions are made up of \emph{almost-invariant sets} \cite{DellnitzJunge99,Froyland05} and are indicative of obstructions to global mixing in the state space \cite{DellnitzEtAl00}. Roughly speaking, approximate level sets of the corresponding eigenvectors (here, left eigenvectors of $P_0$) indicate the location of almost-invariant sets in phase space.
The modern approach to identifying several almost-invariant sets from several leading eigenvectors is to use sparse eigenbasis approximation (SEBA) \cite{Froyland19SEBA} to automatically separate the various approximate level sets. 

The arguments of complex eigenvalues of magnitude close to 1 carry information about the period of slowly decorrelating cycles in dynamical systems \cite{DellnitzJunge99, FroylandEtAl14b, FroylandEtAl21}.
The corresponding complex eigenvectors have high magnitudes on the parts of phase space that support the cycle with period corresponding to the associated complex eigenvalue \cite{FroylandEtAl24,CastroFroyland25}.
While linear response theory shows that the inferred cycle periods are extremely robust to perturbations of the transfer operator arising from isotropic noise \cite{CastroFroyland25}, prior work on low-dimensional fluid flow has shown that small perturbations can have dramatic impact on mixing rates for steady \cite{froyland2017optimal} and periodically forced \cite{froyland2020computation} flows.

\subsection{Eigenvalue optimisation}
\label{sec:evalopt}

We will assume throughout that the matrix $P_0$ is irreducible and aperiodic as a transition matrix of a Markov chain;  this ensures that the eigenvalue 1 is simple and is the only eigenvalue of magnitude 1.
The Gaussian kernel construction from~\eqref{eq:kernel} and~\eqref{eq:markov_mat} guarantees that this property holds for every $\epsilon>0$ if $M$ is connected. Let us denote by $\lambda_0 \in\mathbb{C}$ a simple eigenvalue of $P_0$ such that $|\lambda_0|\lesssim 1$. The corresponding left and right eigenvectors shall be denoted by $u_0,v_0 \in \mathbb{C}^{\tilde{N}}$.
 By matrix perturbation theory \cite{Kato95}, for sufficiently small $\bar\delta > 0$, there is a differentiable family of simple eigenvalues $\lambda_\delta$ of $P_\delta=P_0+\delta\cdot \dot P$  for $\delta\in[0,\bar\delta)$ with $\lim_{\delta\to 0}\lambda_\delta=\lambda_0$. 
Let $u_\delta,v_\delta\in\mathbb{C}^N$  denote the left and right eigenvectors of $P_\delta$ respectively, scaled so that $v_\delta^*v_\delta=1$ and $u_\delta^*v_\delta=1$, where $v_\delta^*$ and $u_\delta^*$ are complex-conjugate transposes of $v_\delta$ and $u_\delta$, respectively.
Standard arguments (e.g.\ Theorem 6.3.12 \cite{HornJohnson13}) show that the derivative of $\lambda_\delta$ at $\delta=0$ is
\begin{equation}
    \label{eq:lamdot} \dot\lambda:=\left.\frac{d\lambda_\delta}{d\delta}\right|_{\delta=0} = u_0^*\dot Pv_0.
\end{equation}
We remark that while we present the framework for matrices in sections \ref{sec:evalopt}--\ref{MagnitudePerturbationSection}, all of the theory extends to the case where the stochastic matrix $P$ is replaced by a Markov Hilbert--Schmidt operator \cite{AntownEtAl22}.

\subsection{Optimal perturbations to increase the eigenvalue magnitude}
\label{MagnitudePerturbationSection}

If we are able to perturb $P_0$ so that the magnitude of $\lambda_0$ increases, this will mean that either (i) the almost-invariant sets associated with $\lambda_0$ will become more invariant ($\lambda_0$ real and close to 1), or (ii) the cycle associated to $\lambda_0$ will decorrelate more slowly ($\lambda_0$ complex with magnitude near 1).
Increasing the magnitude of $\lambda_0$ is equivalent to increasing the real part of $\log\lambda_0$, which we denote as $\Re(\log\lambda_0)$.
This is because $\left.\frac{d|\lambda_\delta|}{d\delta}\right|_{\delta=0}=|\lambda_0|\left.\frac{d(\Re(\log\lambda_\delta))}{d\delta}\right|_{\delta=0}$. Thus, we are concerned with $d(\Re(\log\lambda_\delta))/d\delta = \Re(\lambda_0^{-1}\dot\lambda)=\Re(\bar\lambda_0\dot\lambda/|\lambda_0|^2)$, which using \eqref{eq:lamdot}, becomes
\begin{equation}
    \label{eq:magderiv}
\left.\frac{d(\Re(\log\lambda_\delta))}{d\delta}\right|_{\delta=0}=\frac{\Re\left(\bar\lambda_0\cdot u_0^*\dot P v_0\right)}{|\lambda_0^2|}.
\end{equation}

Following Section 5 of \cite{AntownFroyland18}, the main optimisation problem we therefore wish to solve is
\begin{eqnarray}
    \label{eq:optprob}
\max_{\dot P\in \mathbb{R}^{\tilde{N}\times \tilde{N}}}&& \Re\left(\bar\lambda_0\cdot u_0^*\dot P v_0\right)\\
\label{eq:zerosum}\mbox{subject to}&& \dot P \mathbf{1}=\mathbf{0}\\
\label{eq:norm1}&&\|\dot P\|_F = 1\\
\label{eq:boundary}&&\dot P_{ij}=0 \mbox{ if $(i,j)\in D$},
\end{eqnarray}
for $i,j = 1,\dots,\tilde{N}$, where $D \subseteq \{1, \dots, \tilde{N}\}^2$ is a subset of disallowed transition pairs.
Equation \eqref{eq:zerosum} ensures $\dot P$ has zero row sums, equation \eqref{eq:norm1} normalises the perturbation matrix with the Frobenius norm so that we optimise the infinitesimal perturbation $\dot{P}$ in a unit ball in the tangent space of the feasible set of row-stochastic matrices with base point $P_0$. Condition \eqref{eq:boundary} disallows perturbations at entries of $P_0$ that are specified by the index set $D$. 

Denoting 
\begin{equation*}
S_{ij}=\Re\left(\bar\lambda_0 (\bar{u}_0)_i (v_0)_j\right),~~~ i,j = 1,\dots,\tilde{N}
\end{equation*}
and $D_i^c=\{k:(i,k)\notin D\}$, equation (69) in \cite{AntownFroyland18} shows that the optimal $\dot P$ has elements given by 
\begin{equation}
    \label{eq:optmagdotP}
    \dot P_{ij} = \frac{1}{C} \left(S_{ij} - \frac{1}{|D_i^c|}\sum_{k\in D_i^c} S_{ik}\right),~~~ i,j = 1,\dots,\tilde{N},
\end{equation}
when $(i,j) \notin D$, and $\dot{P}_{ij}=0$ otherwise. 
The term $C>0$ in \eqref{eq:optmagdotP} is a scaling factor chosen so that $\|\dot P\|_F=1$.

\subsection{Optimal perturbations to increase the eigenvalue argument}
\label{FrequencyPerturbationSection}

For complex $\lambda_0$ if we perturb $P_0$ so that the argument of $\lambda_0$ increases, we will be increasing the rate at which this cycle rotates with each application of $P_0$.  In other words we will decrease the period of the cycle or increase its frequency.
Denote the the imaginary part of $\log \lambda_0$ by  $\Im(\log\lambda_0)$.
Since 
$\Im(\log\lambda_0)=\arg \lambda_0$, to increase the argument of $\lambda_0\in\mathbb{C}$ we may equivalently increase $\Im(\log\lambda_0)$.
Following the strategy from the previous subsection, we replace \eqref{eq:magderiv} with 
\begin{equation}\label{eq:freqderivative}
\left.\frac{d(\Im(\log\lambda_\delta))}{d\delta}\right|_{\delta=0}=\frac{\Im\left(\bar\lambda_0\cdot u_0^*\dot P v_0\right)}{|\lambda_0^2|}.
\end{equation}
The optimisation problem is the same as \eqref{eq:optprob}--\eqref{eq:boundary}, except that the objective \eqref{eq:optprob} is replaced with
\begin{equation*}
    \max_{\dot P\in \mathbb{R}^{N\times N}} \Im\left(\bar\lambda\cdot u_0^*\dot P v_0\right).
\end{equation*}
The optimal perturbation matrix $\dot P$ is given by \eqref{eq:optmagdotP}, using $S_{ij}=\Im\left(\bar\lambda_0 (\bar{u}_0)_i (v_0)_j\right)$.

\subsection{Inferring phase space perturbations from $\dot{P}$}
\label{sec:phase_space_perturbations}

In this section we investigate the effects of the perturbation $\dot P$ on the dynamics $\Phi$ through the observation map $F$.
At the outset, we highlight that the following constructions are valid for any perturbation matrix $\dot P$, not only those arising from optimal spectral manipulations.
Before developing the drift for vector-valued observations, suppose that the observation $F:M\to\mathbb{R}$ is scalar-valued and that we do not perform Takens embedding. We observe a length-$(N+1)$ trajectory $\omega_{0},\omega_{1},\ldots,\omega_N\in M$ with $F$ to obtain the scalar observations $x_{0},x_{1},\ldots,x_N\in\mathbb{R}$,  the first $N$ of which are placed in a column vector $\mathbf{x}=[x_0,x_1,\ldots,x_{N-1}]\in\mathbb{R}^{N}$.
Given a $\delta\ge 0$, the vector $P_\delta\mathbf{x}$ is the forward evolution of the scalar observations in $\mathbf{x}$ under one time step of the Markov chain $P_\delta \in \mathbb{R}^{N \times N}$, which corresponds to a time $t=\Delta t$. Similarly, we may interpret the limit of the quotient
\begin{equation*}
\lim_{\delta\to 0}\frac{P_\delta\mathbf{x}-P_0\mathbf{x}}{\delta}=\dot P\mathbf{x}
\end{equation*}
as the derivative with respect to $\delta$ of the $\Delta t$ forward evolution of the vector $\mathbf{x}$.
More precisely $(\dot P\mathbf{x})_i$ represents the infinitesimal $\delta$-response of observation $x_i$ under forward evolution by $\Delta t$.
In this way we may calculate the response of the observations from any perturbation matrix $\dot P$, and in particular the optimal-response of the observations from the optimal $\dot P$.

If $F:M\to\mathbb{R}^d$ is vector-valued, we proceed in the same way for each vector coordinate. This includes the scenario in which Takens delay-embedding is employed, where $F$ is replaced by $F_{Q,L}:M\to\mathbb{R}^{\tilde{d}}$, which is vector-valued for $Q \neq 1$.
Each $F_{Q,L}(\omega_i)=\tilde{x}_i$ lies in $\mathbb{R}^{\tilde{d}}$, and we form an $\tilde{N}\times \tilde{d}$ matrix 
$$\tilde{X}:=\begin{bmatrix}
    \tilde{x}_0\\
    \tilde{x}_2\\
    \vdots\\
    \tilde{x}_{\tilde{N} -1}
\end{bmatrix}$$ with observations forming rows.
To compute the derivative with respect to $\delta$ of the $\Delta t$ forward evolution of the vector-valued observations \emph{along each coordinate separately}, we compute $\dot P\tilde{X}\in\mathbb{R}^{\tilde{N}\times \tilde{d}}$.
In particular, writing 
\begin{equation}
\label{eq:optimal_pert_field}
\begin{bmatrix}
    \dot{\tilde{x}}_0 \\
    \dot{\tilde{x}}_2 \\
    \vdots\\
    \dot{\tilde{x}}_{\tilde{N}-1},
    \end{bmatrix}:=\dot P \tilde{X},
\end{equation}
where $\dot{\tilde{x}}_i$ is a row vector in $\mathbb{R}^{\tilde{d}}$, we obtain the infinitesimal $\delta$-responses of the $\tilde{d}$-vector observations $\tilde{x}$ at all times $i=0,\ldots,\tilde{N}-1$.
The vectors $\{\dot{\tilde{x}}_i\}_{i=0}^{\tilde{N}-1}$ define a discretely sampled vector field, sampled at the observation points $\{\tilde{x}_i\}_{i=0}^{\tilde{N}-1}$. 
When these vectors arise from an optimal perturbation $\dot P$, we call the collection $\{(x_i,\dot{\tilde{x}}_i)\}_{i=0}^{\tilde{N}-1}$, an \textit{optimal-response vector field} of the observable $F$.
This optimal-response vector field realises the optimal perturbation matrix $\dot P$ in terms of the original observations from which $P_0$ was built, and allows us to make visualisations that enable dynamic interpretations of the optimal observation responses. Furthermore, the matrix $\dot P$ produced by the optimisation carried out with the observation function $F$ can be then be mapped onto a different set of observations given by a function $F':M \rightarrow \mathbb{R}^{d'}$ of interest, by replacing $\tilde{X}$ in \eqref{eq:optimal_pert_field} with observations from $F'$.  This will be done in section \ref{ENSOPerturb}, where optimisations carried out with sea-surface temperature can be mapped to optimal observation responses for the wind field. 
We note that a formula related to~\eqref{eq:optimal_pert_field} was used in \cite{SantosGutierezEtAl25} to reconstruct perturbation vector fields in the context of Ulam-type methods---see eq.~(43) in that reference which can be interpreted as a special case of~\eqref{eq:optimal_pert_field} with observations obtained from indicator functions associated with a state space partition.

\subsection{Illustration of optimal frequency-increasing perturbation for the circle rotation}

We now return to analysing the irrational rotation on the circle introduced in Example~\ref{sect:circlerotation}. 
We compute an optimal infinitesimal  perturbation for this system to increase the speed of the rotation.
Let $P_0$ be the transition matrix constructed for this system as in Example 1, and $\lambda_0 = \lambda^{(+1)}$ be its second eigenvalue with corresponding left and right eigenvectors $u_0$ and $v_0$ respectively. 

Using \eqref{eq:optmagdotP} with $S$ replaced as indicated at the end of section \ref{FrequencyPerturbationSection}, we construct the matrix $\dot{P}$ that maximises the derivative of $\arg{\lambda_\delta}$ at $\delta=0$, using the disallowed index set $D=D_\tau := \{(i,j):\tau<P_{ij}<1-\tau\}$.
This set disallows perturbations of transitions with conditional probabilities close to 0 or 1. For the irrational rotation of the circle we use a tolerance of $\tau = 10^{-7}$, as was used in \cite{AntownFroyland18}, resulting in a perturbation matrix $\dot{P}$ with elements that are $47\%$ nonzero. For this $\dot{P}$, the derivative for $\arg\lambda_\delta$ at $\delta=0$, given by \eqref{eq:freqderivative}, is 0.4947.
The largest value for which $P_\delta$ is still stochastic is $\delta = 4.810\times 10^{-5}$, and with this value the eigenvalue argument is increased by $2.380\times 10^{-5}$, corresponding to a frequency increase of $3.787\times 10^{-5}$.

\begin{figure}
    \centering
    \includegraphics[height=5.5cm]{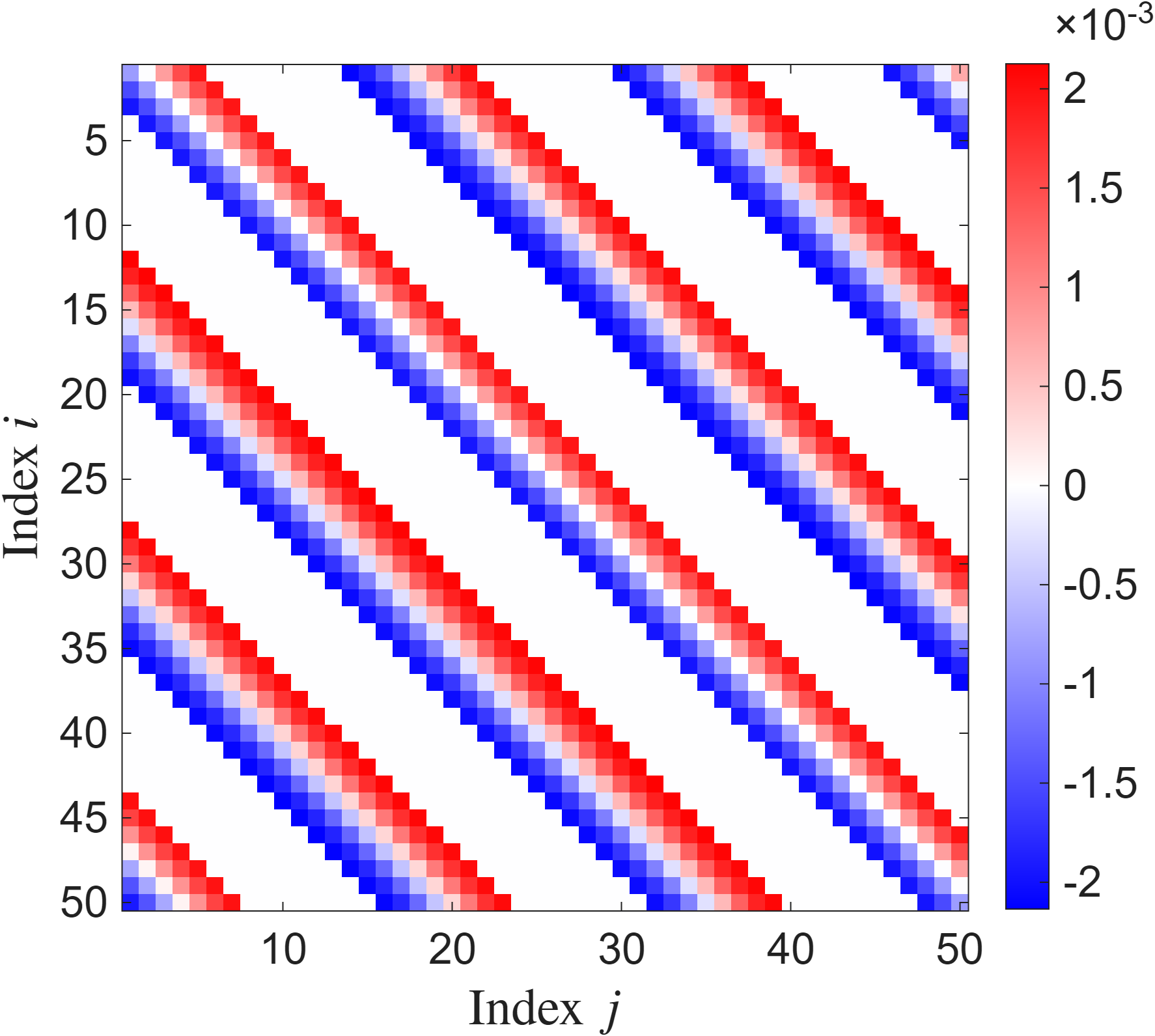}
    \includegraphics[height=5.5cm]{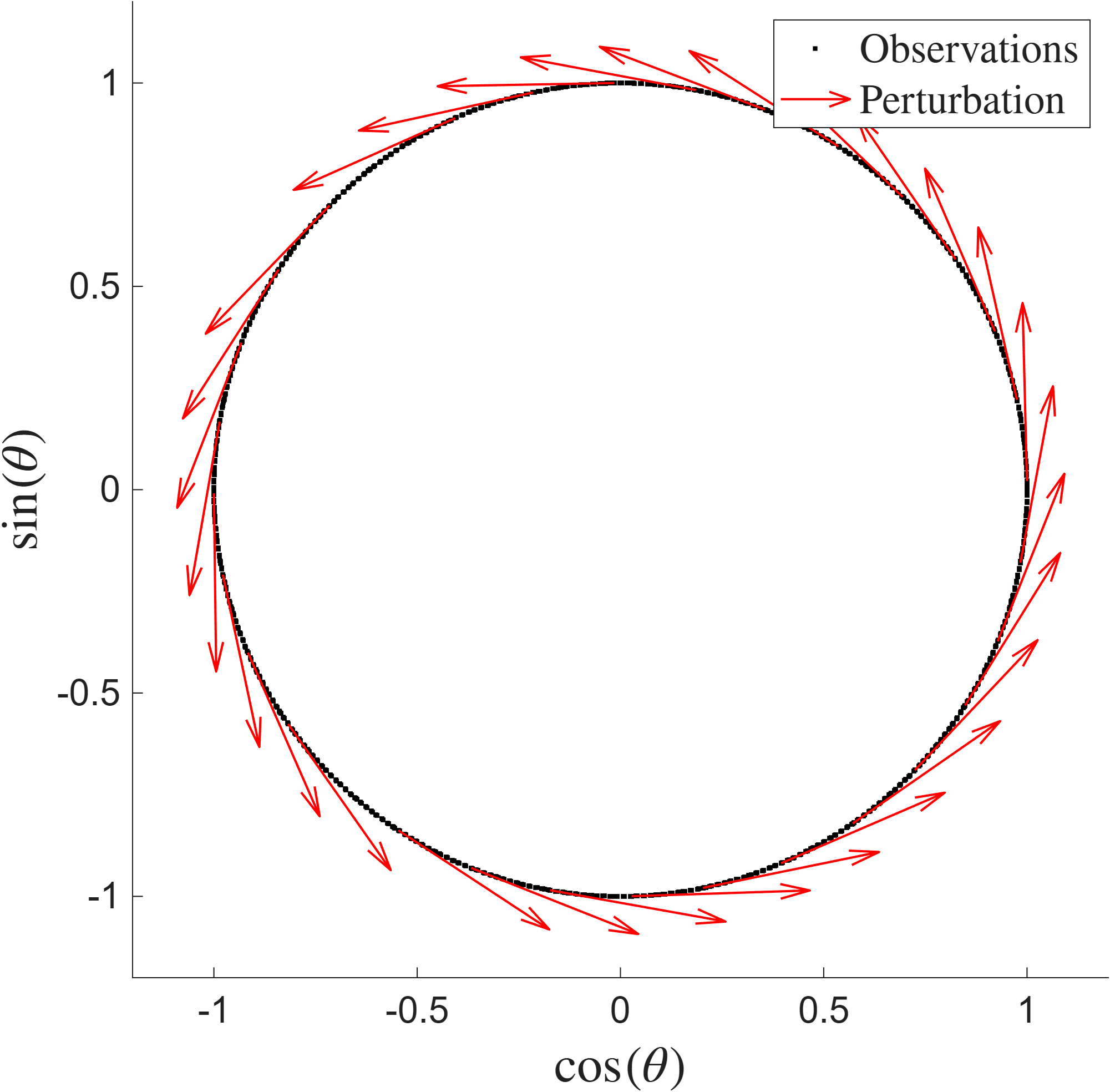}
    \caption{\label{fig:circle_freq_pert}\emph{Left:} Heat map of a $50 \times 50$ diagonal block of the $1000 \times 1000$ perturbation matrix $\dot{P}$. 
    \emph{Right:}  The optimal-response vector field (red arrows), for the frequency-increasing perturbation, originating on the corresponding observational data points. The circle rotation direction is anticlockwise.}
\end{figure}

Figure~\ref{fig:circle_freq_pert}(left) shows that along each diagonal band of $\dot{P}$, positive perturbations are made in the top half of the band and negative in the bottom half of the band.
This promotes motion in the forward direction, increasing the speed. 

Recall in Example~\ref{sect:circlerotation} that the transition matrix $P_0$ was constructed from the observations $x_i = F(\theta_i), i = 0,\dots,N$ using the map $F(\theta) = [\cos \theta, \sin \theta]$. To visualise the effects of the perturbation in the observation space $\mathbb{R}^2$ associated with $F$, we compute the matrix $\dot{P}X \in \mathbb{R}^{N\times2}$ whose rows are the optimal-response 2-vectors $\{\dot{x}_i\}_{i=0}^{N-1}$ as defined in \eqref{eq:optimal_pert_field}.
The optimal-response vector field $\{x_i,\dot{x}_i\}_{i=0}^{N-1}$ is shown in Figure~\ref{fig:circle_freq_pert}(right), where we represent the observations as points and the optimal-response vectors as arrows.
We recover a vector field pointing tangentially to the circle and in the direction of the rotational motion. 
Note that the rotation dynamics occurs on the geometric circle in $\mathbb{R}^2$.
We emphasise that there is no \emph{a priori} constraint on the optimal-response vector field that requires it to be tangential to the geometric circle; {this is a direct outcome of the optimisation of $\dot P$ and its vector-field representation in observation space.}

We now begin a series of case studies, showcasing the power and utility of our optimal-response vector field constructions.

\section{Case study 1: A fast-slow circle rotation}
\label{sec:fast_slow_circle}

\begin{figure}
    \includegraphics[height=4.45cm]{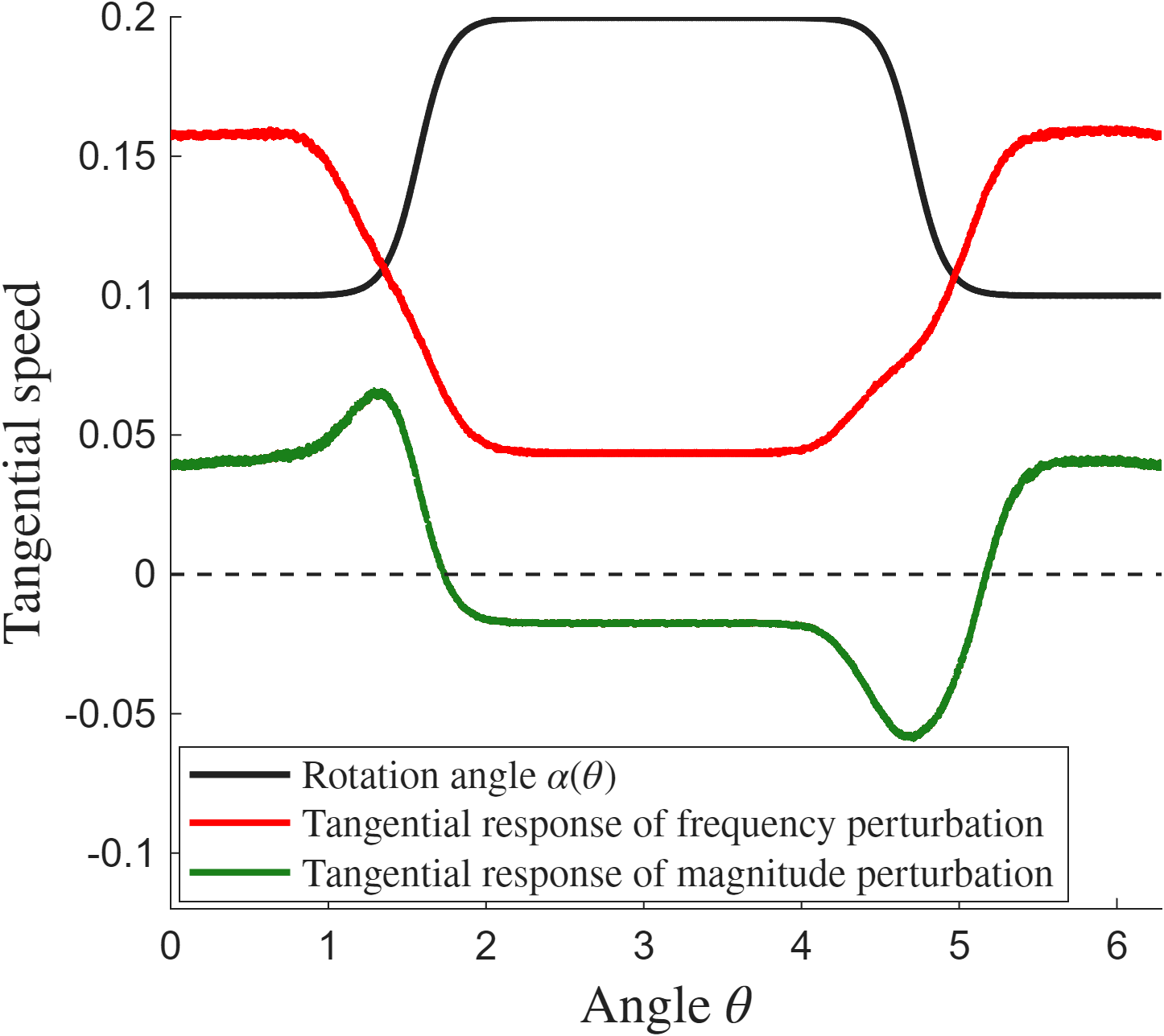}
    \includegraphics[height=4.45cm]{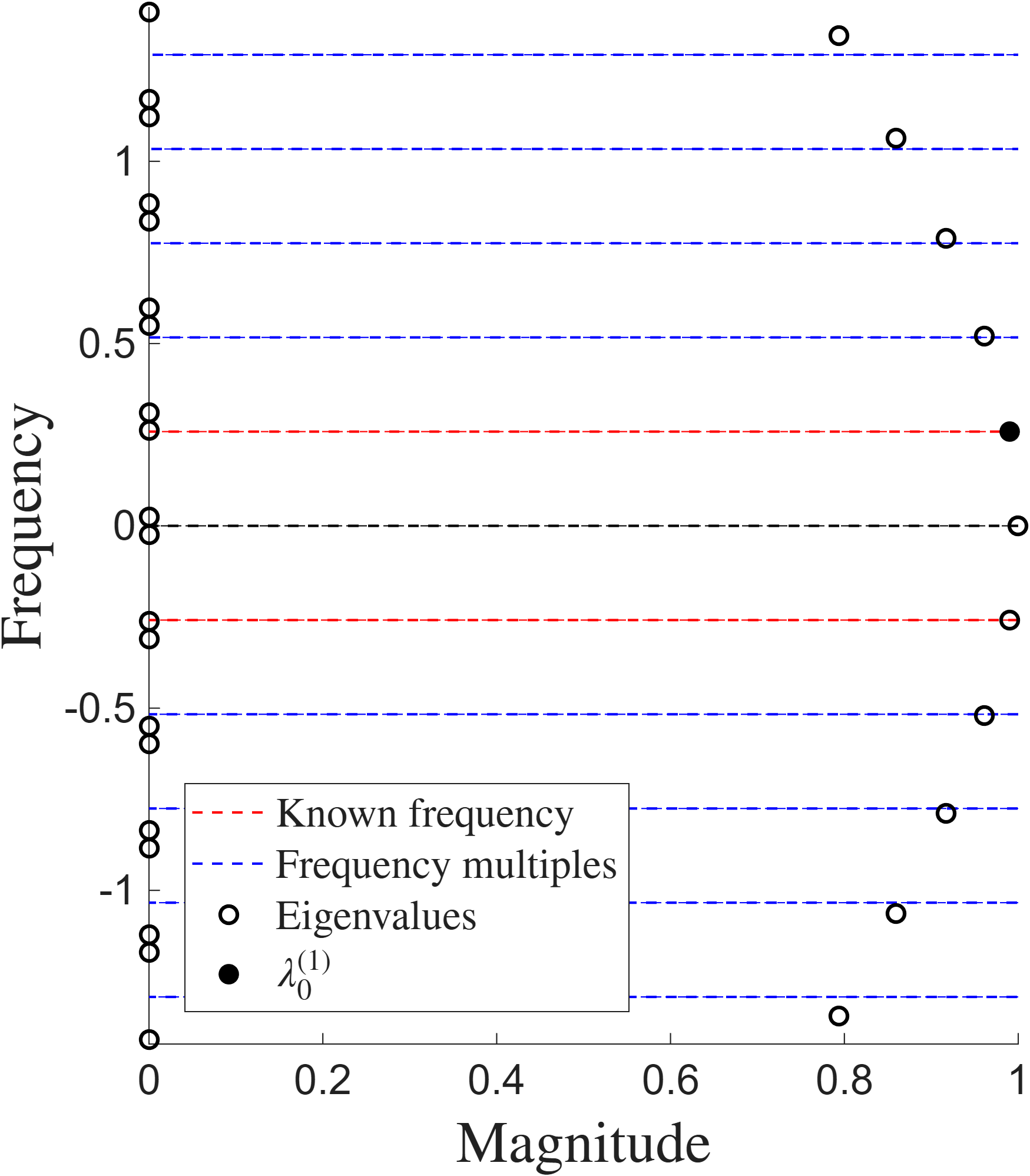}
    \includegraphics[height=4.45cm]{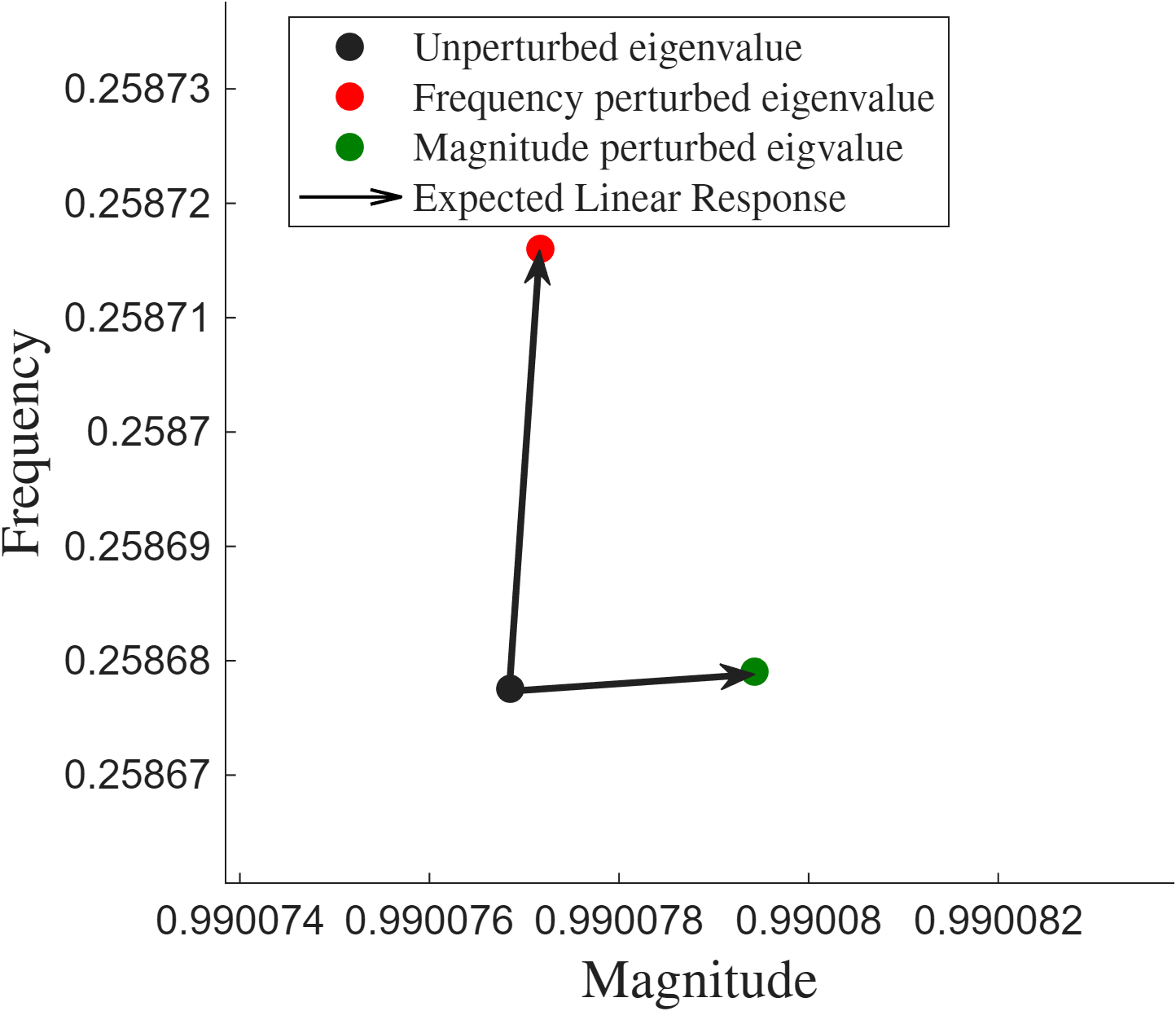}
    \caption{
    \emph{Left:} The variable rotation angle $\alpha(\theta)$ (black line) plotted as a function of the angle $\theta$. Also plotted are the magnitudes of the orthogonal projections of the optimal-response vector fields onto the tangent spaces of the instantaneous velocities in $\mathbb{R}^2$, for the frequency-increasing perturbation (red) and magnitude-increasing perturbation (green). 
    \emph{Centre:} Magnitudes and frequencies of the eigenvalues of $P_0$ constructed using diffusion parameter $\epsilon = 0.1$. \emph{Right:} The unperturbed eigenvalue $\lambda_0^{(1)}$ (black) compared to the perturbed eigenvalues $\lambda_\delta^{(1)}$ associated with the frequency-increasing perturbation (red) and magnitude-increasing perturbation (green).}
    \label{nonuniform}
\end{figure}

Consider the circle rotation $\Phi^{\Delta t}: S^1 \rightarrow S^1$ defined by the map $\Phi^{\Delta t}(\theta) = \theta + \alpha(\theta) \Delta t \mod 2\pi$, where
\begin{equation*}
\alpha(\theta) = 
    \alpha_1 + \frac{1}{2}(\alpha_2-\alpha_1)\left(1+\tanh(5(\theta-\pi/2))\mathbf{1}_{\{0\leq\theta<\pi\}}(\theta)
    -\tanh(5(\theta-3\pi/2))\mathbf{1}_{\{\pi\leq\theta<2\pi\}}(\theta)\right).
\end{equation*}
The constants $\alpha_1 = 0.1$ and $\alpha_2 = 0.2$ are rotation speeds, and $\Delta t = \frac{1}{12}$ is the sampling interval. This map describes an irrational rotation of the circle, where the rotation angle is $\alpha_1$ in most of the left half circle $\theta \in [\frac{\pi}{2},\frac{3\pi}{2})$ and is $\alpha_2$ in most of the right half circle. The hyperbolic tangent function introduces a smooth transition between the two speeds. This is shown in Figure \ref{nonuniform}(left), by plotting $\alpha(\theta)$ as a function of $\theta$ (black curve).

The period of this rotation is
\begin{equation*}
    T = \Delta t \int_{S^1} \frac{1}{\alpha(\theta)} d \theta
    = \Delta t  \cdot (15 \pi - \ln(2))\\
    \approx 3.8692 \text{ time units,
    }
\end{equation*}
which corresponds to a frequency of $\nu = \frac{1}{T} \approx 0.2584$.

The system is observed by $F:S^1 \rightarrow \mathbb{R}^2$, defined by $F(\theta) = [\cos(\theta),\sin(\theta)]$. We collect $N+1 = 4000$ observations $x_i = F(\Phi^{i\Delta t}(\theta_0)), n = 0,\dots, N$ of the trajectory initialised at $\theta_0 = 0$. Delay-embedding is unnecessary for this system so we set $Q=1$ and proceed with $\tilde{x}_i = x_i$ as our data for constructing the Markov transition matrix. For the construction of $P_0$, the diffusion parameter is set to $\epsilon = 0.1.$

Figure \ref{nonuniform}(centre) shows how the frequency $\nu^{(+1)} = 0.2586$, associated with the leading nontrivial eigenvalue $\lambda^{(+1)}$ of $P_0$, is approximately the theoretical frequency $\nu=0.2584$. Successive eigenvalues have frequencies which are approximately integer multiples of $\nu$.

\subsection{Perturbations of the non-uniform circle rotation}

We investigate the effects of increasing (i) the frequency and (ii) the strength of the non-uniform rotation cycle. This is done by perturbing $P_0$ to increase the frequency and magnitude of $\lambda_0=\lambda^{(+1)}$, respectively. For both perturbations, $\dot P$ is constructed with a tolerance of $\tau = 10^{-7}$, resulting in a perturbation matrix with $22\%$ nonzero elements. To optimally increase the frequency, we construct the perturbation matrix $\dot P$ which maximises the derivative of $\arg \lambda_\delta$ with respect to $\delta$.
Then according to \eqref{eq:freqderivative}, the derivative of $\arg\lambda_\delta$ at $\delta=0$ is 0.1052.

\begin{figure}
\includegraphics[width=0.49\textwidth]{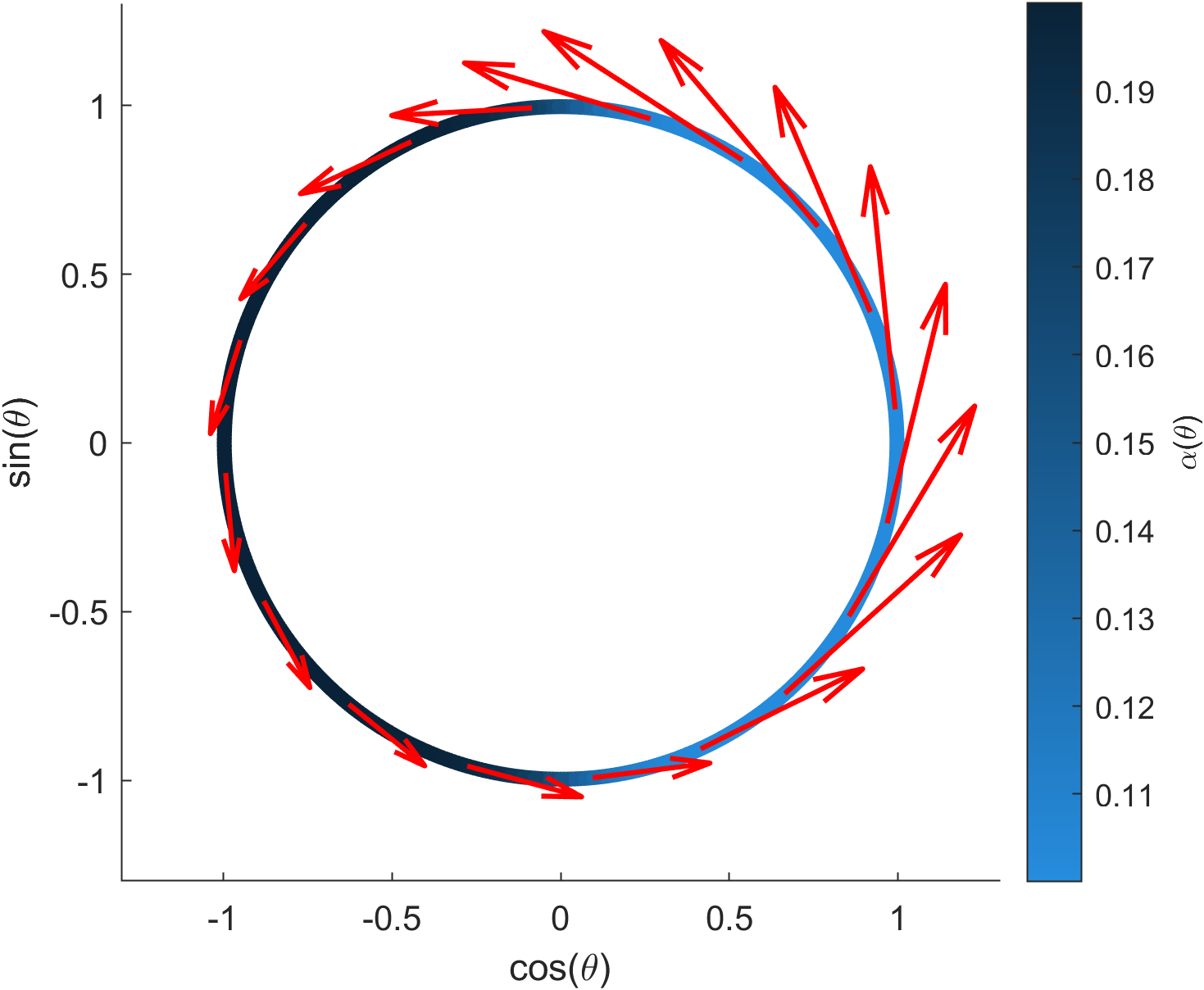} \includegraphics[width=0.49\textwidth]{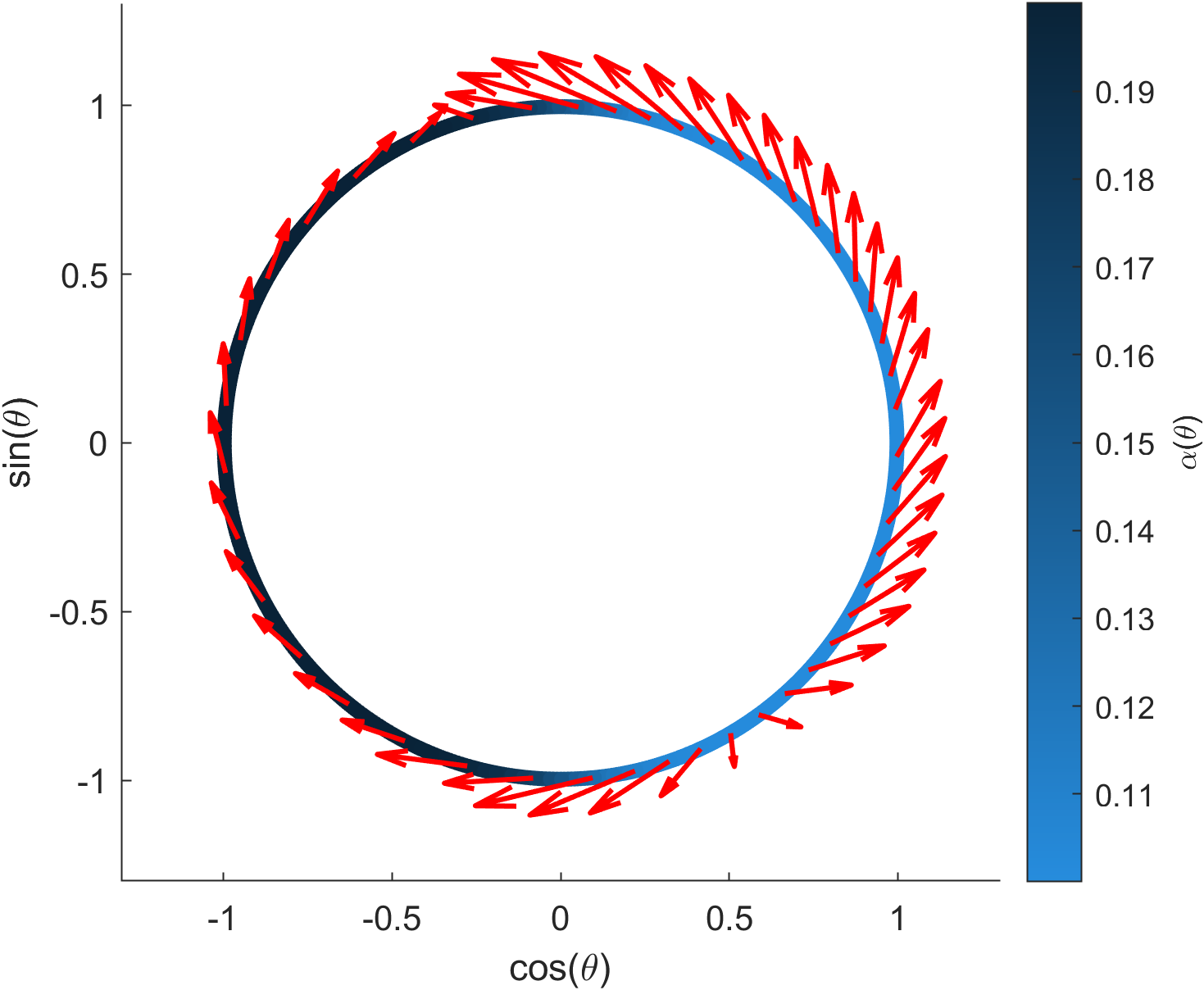}
\caption{\emph{Left:} The observational data is coloured according to the forward difference $\alpha(\theta)$. The arrows represent the optimal-response vector field of the frequency-increasing perturbation at a selection of points on the circle. \emph{Right:} The same information for the magnitude-increasing perturbation. The scaling factor and number of arrows were chosen independently to aid in visualisation, and are not the same between the two figures.}
\label{fig:NonuinformDrifts}
\end{figure}

The corresponding optimal-response vector field is shown in Figure \ref{fig:NonuinformDrifts}(left) in the observation space. 
Similar to the optimal-response vector field for the uniform circle rotation shown in Figure \ref{fig:circle_freq_pert}(right), the vectors point in the direction of the unperturbed rotation and they are close to tangential to the circle. However, in this nonuniform rotation case the optimal perturbation has greater strength during the slow portion of the nonuniform rotation. 
In Figure \ref{nonuniform}(left) in red we  project the vectors in Figure \ref{fig:NonuinformDrifts}(left) onto the tangent space of the circle to compare with the tangential velocity. In Figure \ref{nonuniform}(left), the speed $\alpha(\theta)$ and the optimal perturbation are a function of angle $\theta$.
The overall perturbation is positive, corresponding to positive advancement of the rotation, and the perturbation is larger where $\alpha(\theta)$ is smaller.

This may be explained as follows.
If at each $\theta\in [0,2\pi)$ we perturb the velocity $\alpha(\theta)$ by an amount $\Delta\alpha(\theta)$, Jensen's inequality says
$$\int_0^{2\pi} \frac{1}{\alpha(\theta)+\Delta\alpha(\theta)}\ d\theta \ge \frac{2\pi}{\frac{1}{2\pi}\int_0^{2\pi} \alpha(\theta)+\Delta\alpha(\theta)\ d\theta}.$$
The LHS above is the time taken for the perturbed flow to make one revolution, while the RHS above is the time taken to traverse the circle at the \emph{average perturbed speed}.
Equality holds only when the speed $\alpha(\theta)+\Delta\alpha(\theta)$ is constant (independent of $\theta$). Thus, to minimise the time taken to make one revolution, convexity through Jensen says that one should choose $\Delta\alpha$ to make $\alpha+\Delta\alpha$ more constant, by speeding up more in slow sections and speeding up less in fast sections. This is precisely what we see in Figure \ref{fig:NonuinformDrifts}(left).

Next, we optimally increase the \emph{strength} of the cycle by using the $\dot P$ which optimally increases $|\lambda_\delta^{(+1)}|$. The resulting derivative given by \eqref{eq:magderiv} is 0.0153, which corresponds to a derivative for $|\lambda_\delta|$ at $\delta=0$ of 0.0152. These two values are similar because $|\lambda_0| \approx 1$. The associated optimal-response vector field is shown in Figure \ref{fig:NonuinformDrifts}(right). 
The vectors are again close to tangential; the main difference between the left and right panels of Figure \ref{fig:NonuinformDrifts} is that the vector field in the fast section of the rotation is retrograde relative to the direction of rotation, so that the perturbed rotation is slowed down in this half of the circle.

A possible explanation for this is as follows. 
Firstly, to increase the magnitude of $\lambda_0$ the system will seek to reduce the total $\epsilon$-noise encountered over a cycle from the construction of $P_0$.
At each time step $\epsilon$-noise is incurred, so to reduce overall $\epsilon$-noise, the system would prefer to ``even out'' the noise accumulation by evening out the rotation speed.
Secondly, as seen in Figure \ref{nonuniform}(right), the perturbed eigenvalue moves almost completely in the magnitude direction without altering the frequency of the cycle.
Therefore, the overall period of the perturbed cycle must remain the same as the unperturbed cycle.
Combining these two aspects, for the same convexity reasons mentioned earlier, the optimisation desires to increase speed more in the slow section of rotation relative to the fast section of rotation.
The difference now is that total period must be maintained, necessitating the negative counter-rotation perturbation seen in Figure \ref{nonuniform}(left) as the green curve.
Similarly to the red curve, the green curve is formed by orthogonal projection onto the tangent space of the circle.

\section{Case study 2: the Lorenz flow} \label{section:lorenz}
Consider the Lorenz system \cite{Lorenz63}
\begin{align*}
    \dot{x} &= \sigma(y-x),\\
    \dot{y} &= x(\rho-z)-y,\\
    \dot{z} &= xy-\beta z,
\end{align*}
with standard parameters $\sigma=10, \rho = 28, \beta = \frac{8}{3}$, for which the system possesses the iconic butterfly-shaped chaotic attractor in Figure \ref{LorenzFreqDrift}(left).
{It has been proven that the Lorenz attractor is the support of the SRB measure, $\mu$, of the system \cite{Tucker99} and the dynamics with respect to $\mu$ is mixing \cite{LuzzattoEtAl05}.}

In this section, we perform three separate optimisations. We first focus on the leading complex eigenvalue, which is associated with the transport of an almost-cyclic set around the wings \cite{FroylandEtAl21}, and identify the perturbations that optimally increase the (i) frequency and  (ii) strength of this cycle. We then study the leading nontrivial real eigenvalue, whose eigenvector partitions the attractor into two almost-invariant sets \cite{DellnitzJunge99, FroylandDellnitz03, FroylandPadberg09}, and determine the perturbation that strengthens the almost-invariance property of these sets.

\subsection{Identification of the dominant long-lived cycle in the Lorenz attractor}

We begin by describing a slight variant of the example described in \cite{FroylandEtAl21, CastroFroyland25}. The input data is a trajectory of length $N+1=32000$ with a sampling interval of $\Delta t = 0.01$, for a total timespan of 320 time units. During this time the trajectory completed 426 laps\footnote{We define a lap as the subinterval of the trajectory in which it starts at a local minimum of the $z$ coordinate, rises and falls as it traverses a wing of the attractor and then arrives the next local $z$ minimum.} around one of the wings of the Lorenz attractor. The mean period taken by the trajectory to complete a lap around either wing is $T_{\rm lap} = 0.751$, with standard deviation 0.0947. The empirical frequency $\nu_{\rm lap}=1/T_{\rm lap} = 1.332$ will be the frequency we later compare our eigenvalues to when extracting the dominant cycle.

\begin{figure}
    \centering    \includegraphics[width=0.5\linewidth]{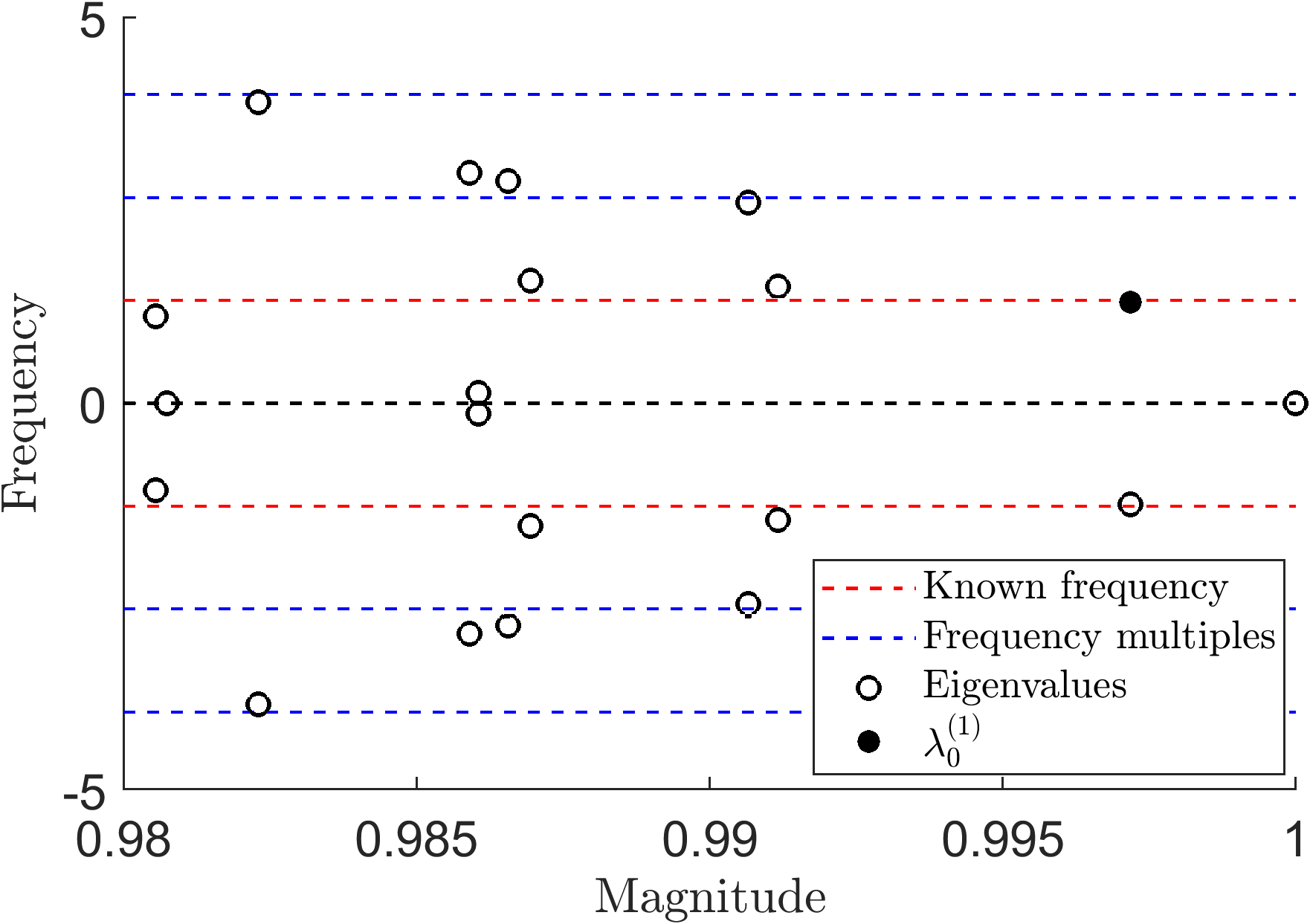}
    \caption{Magnitudes and frequencies of the leading eigenvalues of $P_0$ with $\epsilon = 0.022$, approximating those of the Koopman operator for the Lorenz system.}
    \label{LorenzEigVals}
\end{figure}

Let us denote the state by $\omega = (x,y,z)$. The trajectory is observed under the identity map $F=\mathrm{Id}:\mathbb{R}^3 \rightarrow \mathbb{R}^3$ defined by $F(\omega)=\mathrm{Id}(\omega) = \omega$, without delay-embedding (setting $Q=1$). The Markov matrix $P_0$ was computed from the Lorenz trajectory data using \eqref{eq:markov_mat} with parameter $\epsilon = 0.022$.
As illustrated in Figure \ref{LorenzEigVals}, the leading nontrivial complex eigenvalue is $\lambda_0 = \lambda^{(+1)}=0.997e^{\pm i 0.0822}$, which has a frequency $\nu^{(+1)} = \frac{\arg(\lambda^{(+1)})}{2\pi \Delta t} = 1.309$ within 2\% of the empirical frequency $\nu_{\rm lap} = 1.332$. 
\begin{figure}
    \centering    \includegraphics[width=0.67\linewidth]{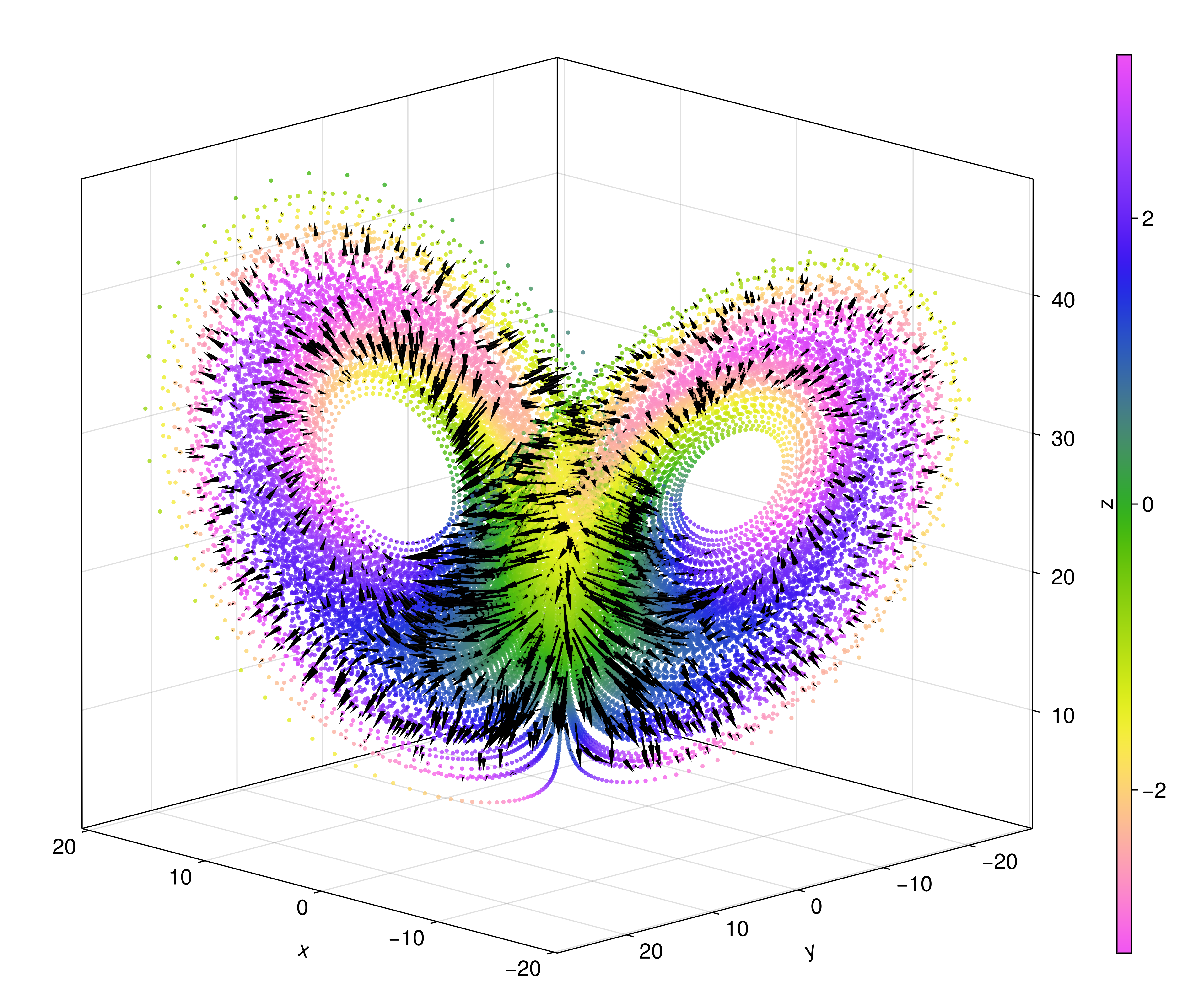}
\includegraphics[width=0.32\linewidth]{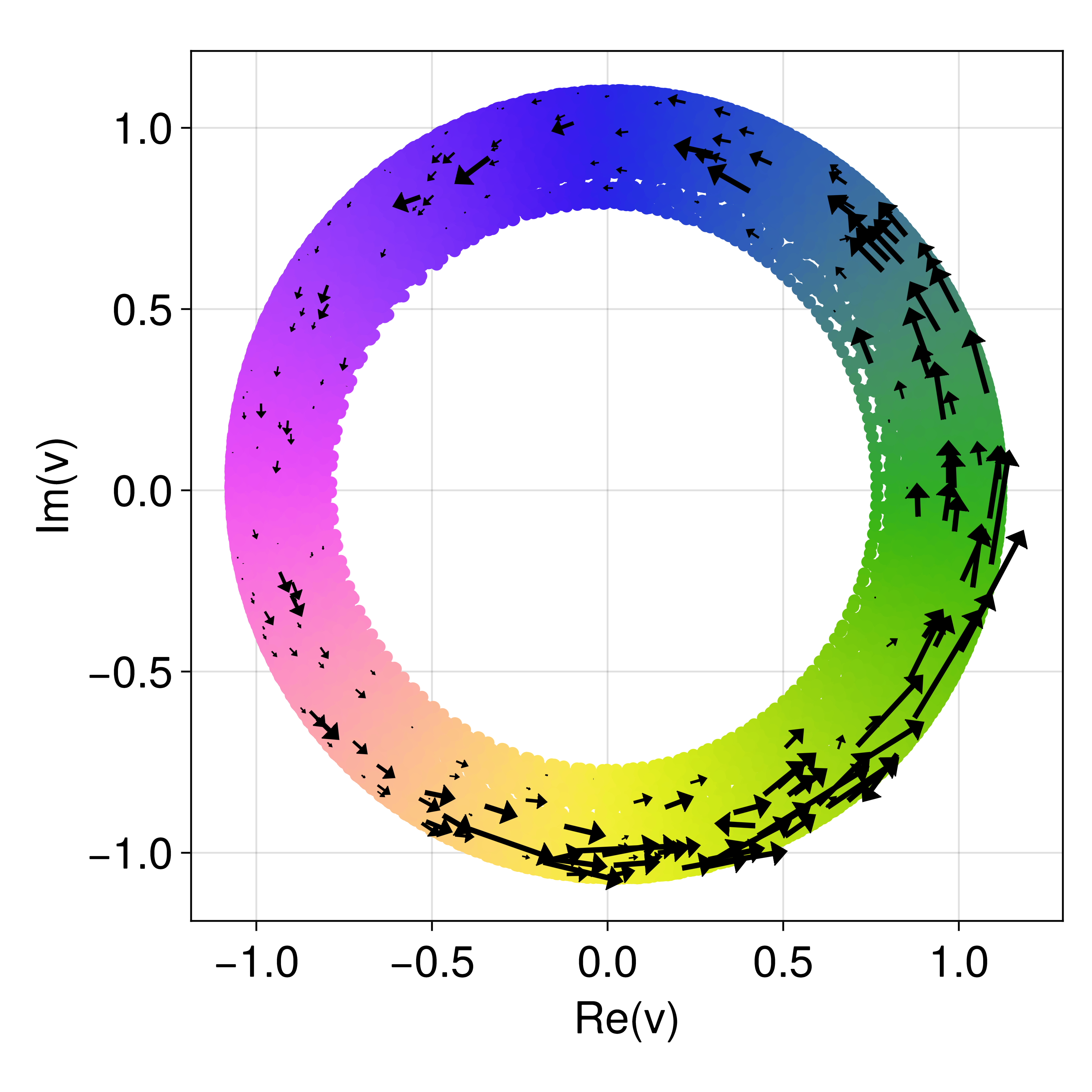}
    \caption{\emph{Left:} The observational data is coloured according to the argument of $v_0$. The arrows represent the optimal-response vector field of the frequency perturbation. For visibility, arrow lengths have been uniformly scaled and only a random sample of arrows are drawn, sampled with weighting inversely proportional to the density of points.  \emph{Right:} Entries of $v_0$ are plotted in the complex plane, coloured by their argument, using the same colour scheme as on the left. The arrows are the pushforward of the vector field on the left, from $\mathbb{R}^3$ to $\mathbb{C}$.
    }
    \label{LorenzFreqDrift}
\end{figure}
We emphasise that the calculation of $\nu_{\rm lap}$ contains a range of frequencies estimated by the standard deviation of 0.15, whereas $\nu$ represents a cycle frequency in a more restricted central part of the wings indicated by the magnitude of the eigenvector;  see \cite{CastroFroyland25} for more details.

\begin{figure}
    \centering    \includegraphics[width=0.67\linewidth]{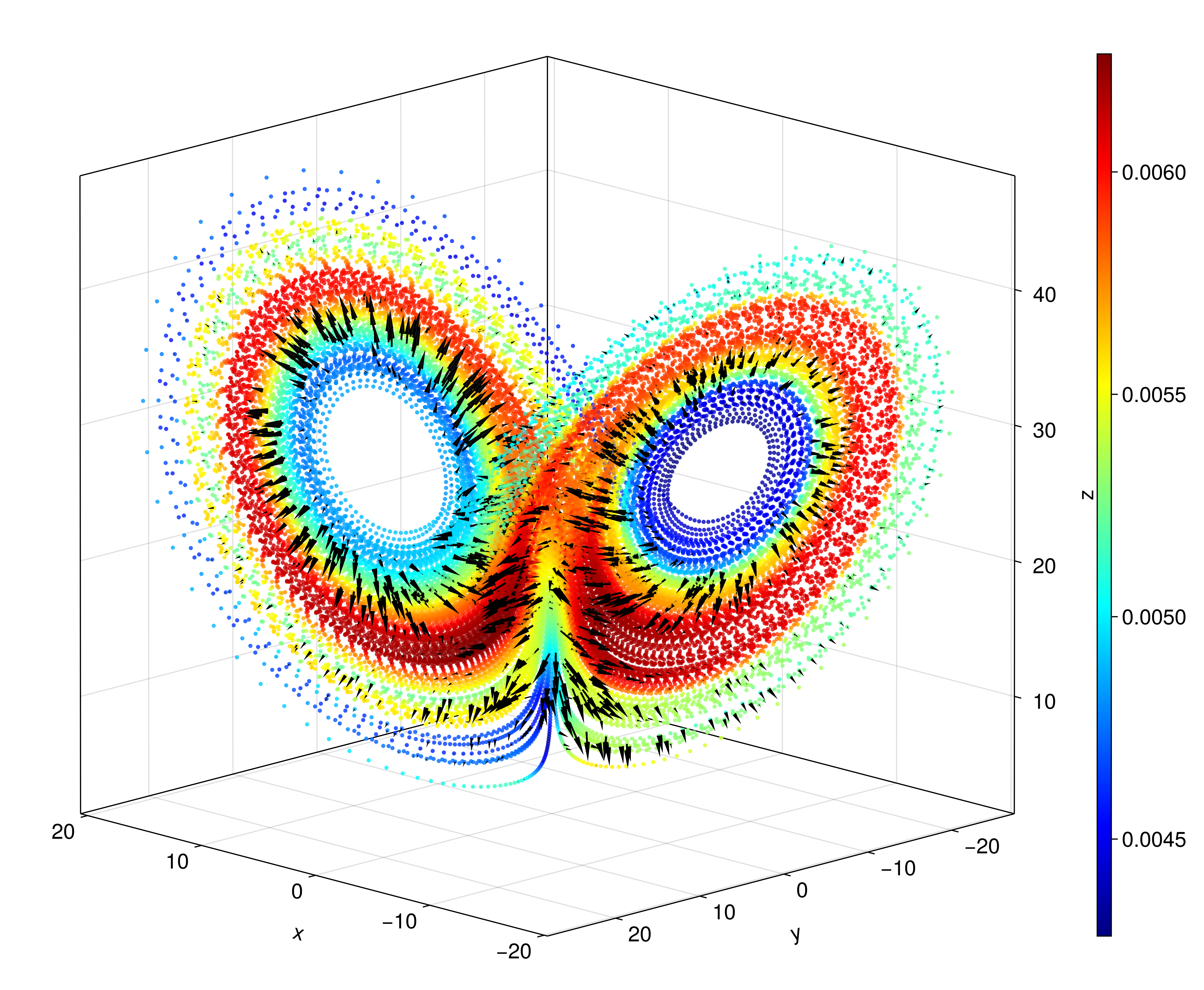}
\includegraphics[width=0.32\linewidth]{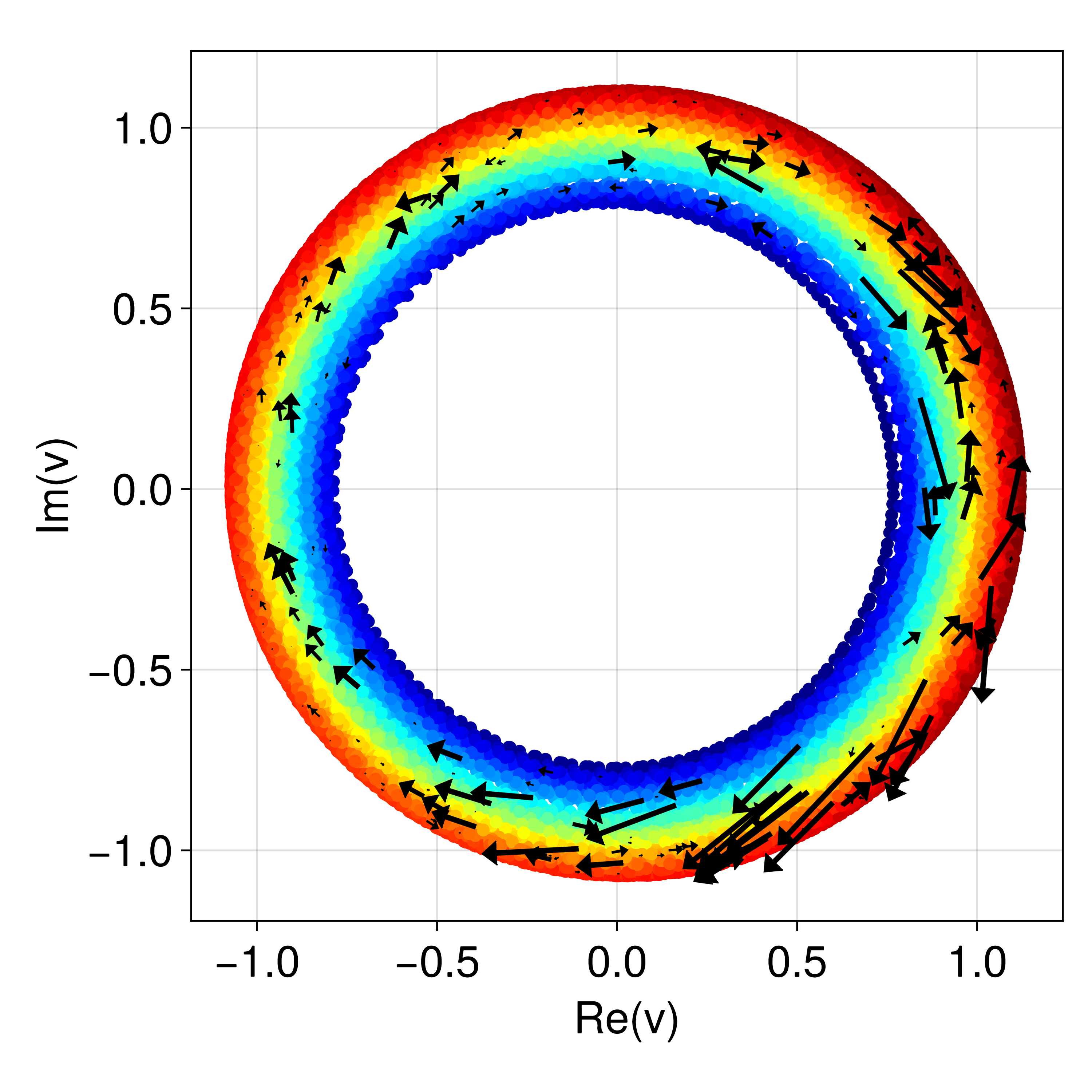}
    \caption{\emph{Left:} The observational data is coloured according to the magnitude of $v_0$. The arrows represent the optimal-response vector field of the magnitude perturbation. Arrow selection and length scaling are as in Figure \ref{LorenzFreqDrift}. \emph{Right:} Entries of $v_0$ are plotted in the complex plane, coloured by their magnitude, using the same colour scheme as on the left. The arrows are the pushforward of the vector field on the left, from $\mathbb{R}^3$ to $\mathbb{C}$.}
    \label{LorenzMagDrift}
\end{figure}

Similar to the results obtained in \cite{FroylandEtAl21}, Figure \ref{LorenzFreqDrift}(right) shows that the leading complex eigenvector forms a circle when its entries are plotted in the complex plane.
This circle rotates anticlockwise with a period $T=1/\nu=0.764$. Similar to the experiments in \cite{FroylandEtAl21,CastroFroyland25}, in Figure \ref{LorenzFreqDrift}(left), the observational data is coloured according to the argument of the eigenvalue at each point. The argument cycles through its colour scale as the trajectory travels around either wing, and by comparing the left and right figures, it can be seen that the trajectory completes a lap on the attractor each time the eigenvector completes a full rotation in the complex plane. Figure \ref{LorenzMagDrift}(left) displays the magnitude of $v_0$, with the peak magnitude occurring on the ``figure-of-8" red central band.
The frequency $\nu$ of the cycle identified by the eigenvalue of $v_0$ is most accurate on this band; these results duplicate the findings in \cite{CastroFroyland25}. 

\subsubsection{Optimally decreasing the period of the cycle around the wings}
\label{sec:l63_freq_pert}

We now wish to increase the cycle frequency.
To achieve this, we construct $\dot P$ using \eqref{eq:optmagdotP}, making the substitution for $S$ at the end of section \ref{FrequencyPerturbationSection}. We use a tolerance of  $\tau = 10^{-10}$ to construct $\dot{P}$. Such a small tolerance is necessary for this system because the small diffusion parameter $\epsilon = 0.022$ means $P_0$ is very sparse (2\% nonzero).
The derivative of $\arg \lambda_\delta$ with respect to $\delta$, computed with \eqref{eq:freqderivative} is 0.0283.

The arrows in Figure \ref{LorenzFreqDrift}(left) represent the optimal-response vector field $\dot P X$ for this frequency-increasing perturbation. One might intuitively expect these arrows to point in the direction of flow, but we instead see them point in the direction of increasing eigenvector argument, 
i.e.\ the arrows are orthogonal to level sets of the argument.

To make it transparent that the optimal vectors contained in $\dot X$ are indeed orthogonal to level sets of the eigenvector argument, we consider the mapping $\Psi_{v_0}:\mathbb{R}^3\to \mathbb{C}$ defined discretely on the trajectory $\{x_i\}$ by $\Psi_{v_0}(x_i)=v_{0,i}$, and smoothly extend $\Psi_{v_0}$ to a open neighbourhood of the Lorenz attractor.
Identifying $\mathbb{C}$ with $\mathbb{R}^2$, we may then estimate $D_{x_i} \Psi_{v_0}:\mathbb{R}^3\to \mathbb{R}^2$ and for each $x_i$, compute $D_{x_i}\Psi_{v_0}(\dot X_i)\in \mathbb{R}^2$, where $\dot X_i\in \mathbb{R}^3$ denotes the optimal vector at $x_i$. Numerically, this process is particularly straightforward using the Julia package \verb"RadialBasisFunctions.jl" \cite{RadialBasisFunctions_jl}.
{When interpreting the results below, it should be kept in mind that a general vector field need not be projectable under a non-invertible map such as $\Psi_{v_0}$.
That is, it is possible that $D_{x_i} \Psi_{v_0}(\dot X_i)$ and $D_{x_j} \Psi_{v_0}(\dot X_j)$ (which are represented by arrows anchored at the points $\Psi_{v_0}(x_i)$ and $\Psi_{v_0}(x_j)$ in $\mathbb R^2$, respectively) are not equal even if $\Psi_{v_0}(x_i) = \Psi_{v_0}(x_j)$.}

Figure \ref{LorenzFreqDrift}(right) displays this pushforward of the optimal-response vector field.
{The smoothness of the numerically projected vector field in Figure~\ref{LorenzFreqDrift}(right) suggests that the optimal-response vector field $\dot X_i$ \emph{is} projectable under $\Psi_{v_0}$, which is a non-trivial property.}
Furthermore, the arrows point tangentially in the rotation direction of the eigenvector (and therefore orthogonally to level sets of the argument), which clearly speeds up the associated cycle. The pushforward response vector field also contains some information about the vector magnitudes, since the long arrows in the intersection between the wings are mapped to long arrows in the yellow-green part of the eigenvector.
These long arrows indicates that the optimal-response field is strongest in this intersection, which can be explained by the fact that this region is where the cycle is slowest. Indeed, Figure \ref{LorenzFreqDrift}(left) shows that in this intersection the colour progresses from pink-yellow to green-blue, which accounts for roughly half of the cycle in Figure \ref{LorenzFreqDrift}(right), despite being only a small piece of the attractor. This slow sector is where the optimal-response vector field is strongest, which is similar to the results we obtained for the frequency-increasing perturbation for the non-uniform circle rotation in Figure \ref{fig:NonuinformDrifts}(left).

\subsubsection{Optimally strengthening the cycle around the wings}

We repeat the previous experiment but instead perturb to optimally increase the strength of the cycle, by computing $\dot{P}$ using \eqref{eq:optmagdotP} to optimally increase the derivative of $\Re(\log(\lambda_0))$ with respect to $\delta$. Again, the tolerance was set to $\tau = 10^{-12}$. The value of the derivative \eqref{eq:freqderivative} is 0.0163, which is approximately the derivative of $|\lambda_\delta|$ with respect to $\delta$, since $|\lambda_0| \approx 1$.

The arrows in Figure \ref{LorenzMagDrift}(left) represent the optimal-response vector field $\dot P X$. 
The arrows fall into two main types. Firstly those in the inner parts of the wings ``inside'' the red band, which point toward the red band. This is to move the trajectories closer to the red band, where the cycle has less variation in period. Secondly, those at the outer edges of the wings, which point further toward the outer wing edges. This is the most efficient way to get back into the red band because the outer wing edges will be reinjected into the inner parts of the wings, which will then move outward toward the red band. 

As for the frequency perturbation, we again push the response vector field forward onto the eigenvector in $\mathbb{C}$, to obtain the arrows in Figure \ref{LorenzMagDrift}(right). 
Unlike Figure \ref{LorenzFreqDrift}(right), which shows a smooth vector field, Figure \ref{LorenzMagDrift}(right) appears to be very rough.
{This is likely due to the optimal-response vector field associated with strengthening the dominant cycle of the Lorenz flow not being projectable under the map $\Psi_{v_0} \colon \mathbb R^3 \to \mathbb C$ associated with the eigenfunction identifying that cycle (in contrast to the frequency perturbation case from section~\ref{sec:l63_freq_pert}).} 

For example, there are two circular curves per wing on which $|v_0|$ takes a given value, one near the inner edge of the attractor and one near the outer edge.
On a single wing, when these two curves meet a chosen level set of $\arg{v_0}$, one obtains two points at which $v_0$ is equal. 
Including both wings, there are thus four points on the attractor that are sent to a single point in $\mathbb{C}$.
For the vector field in Figure \ref{LorenzMagDrift}(right) to be smooth, we require that the vector field on the attractor respect the symmetries: (i) the left to right wing symmetry and (ii) the inner band to outer band symmetry. One can see from Figure \ref{LorenzMagDrift}(left)that this is not the case, and therefore the vector field is not smoothly projectable onto the annulus in $\mathbb{C}$, rendering information from this alternate representation limited.

\subsection{Identifying the almost invariant-sets which govern the switching between attractor wings}

Another interesting feature of the Lorenz attractor are the almost-invariant sets, which arise from the fact that trajectories tend to stay for long durations on a wing before switching to the other wing. Sets that optimise this metastability or almost-invariance property have been studied with transfer operators \cite{DellnitzJunge99, FroylandDellnitz03, FroylandPadberg09} and here we ask the question:  how should the flow be perturbed to maximally strengthen the almost-invariance?

\begin{figure}
    \centering
    \includegraphics[width=0.5\linewidth]{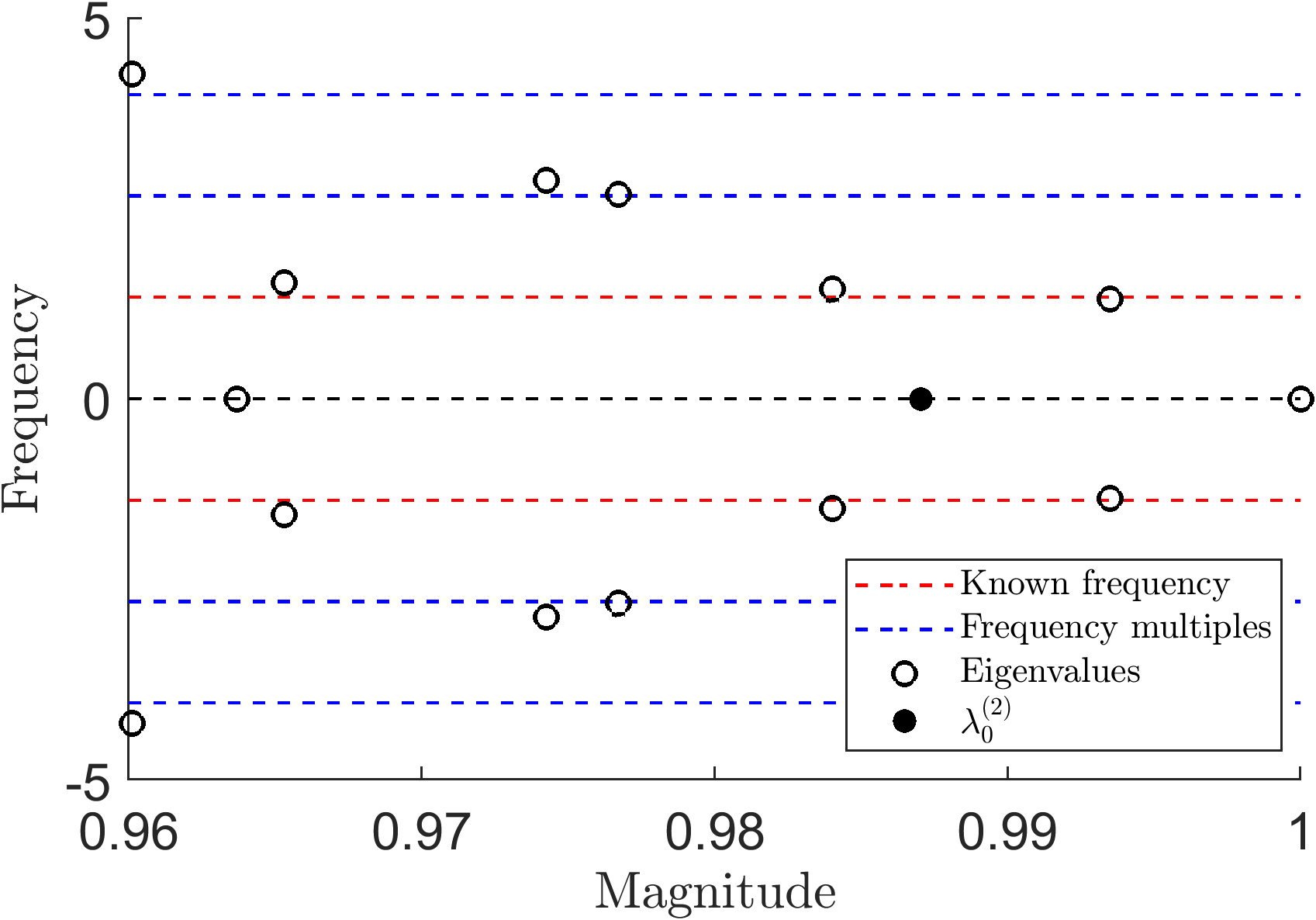}
    \caption{The magnitudes and frequencies of the leading eigenvalues of $P_0$ with $\epsilon = 0.039$, {representing approximate eigenvalues of the} Koopman operator for the Lorenz system.}
    \label{LorenzRealEigvals}
\end{figure}

\begin{figure}
    \centering
    \includegraphics[width=0.67\linewidth]{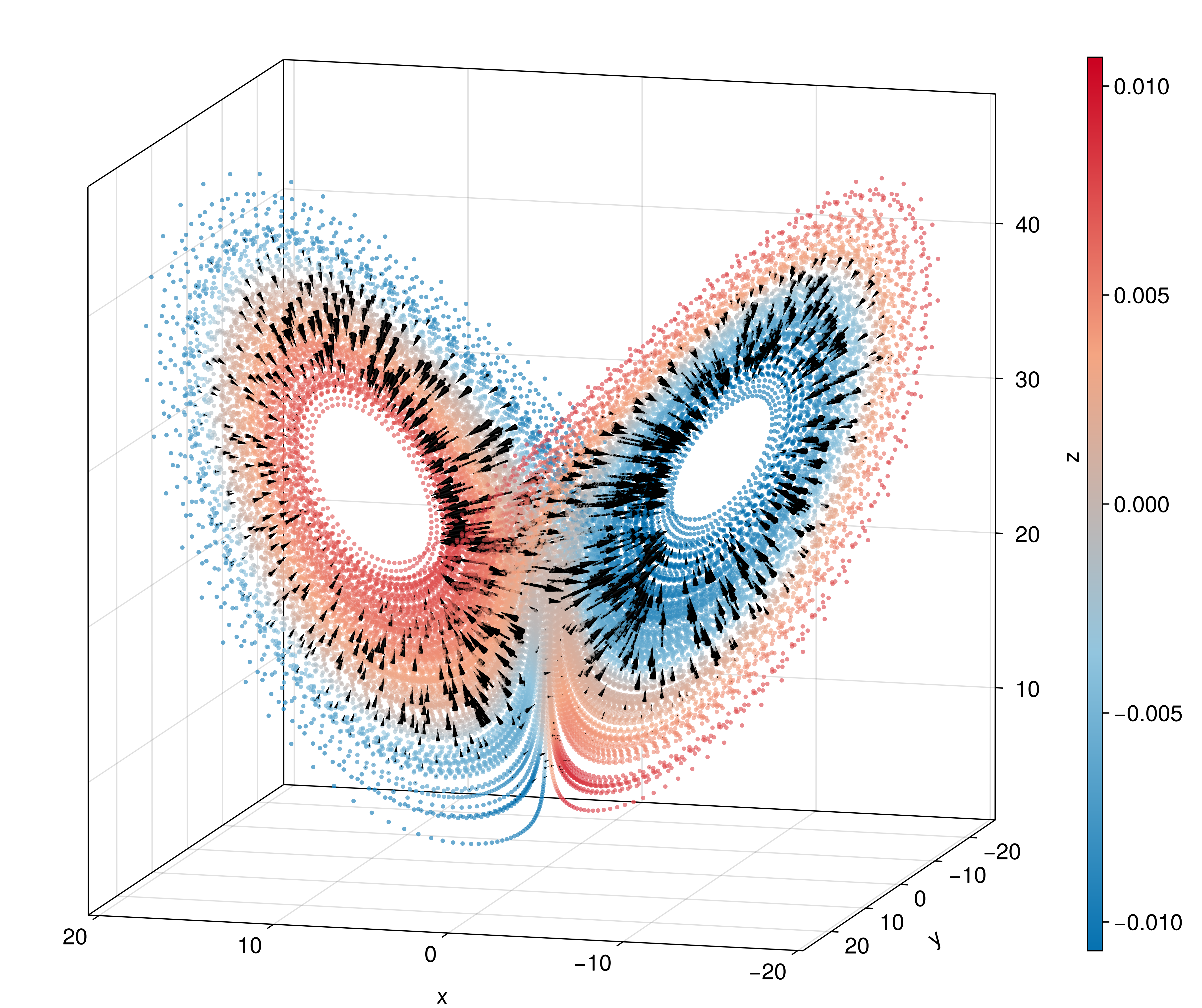}
    \caption{The observational data is coloured according to the eigenvector $v_0$. The arrows represent the optimal-response vector field of the perturbation which decreases the transport between the almost-invariant sets. Arrow selection and length scaling are as in Figure \ref{LorenzFreqDrift}.
    }
    \label{LorenzRealDrift}
\end{figure}

We increase the diffusion parameter to $\epsilon = 0.039$, to obtain a real leading nontrivial eigenvalue (see Figure \ref{LorenzRealEigvals}) and a corresponding  eigenvector $v_0$ that takes positive values on one of the two almost-invariant sets and negative values on these other (see Figure \ref{LorenzRealDrift}).
The dynamics underlying each almost-invariant set is as follows. Points in the inner part of a wing spiral around, gradually increasing in radius, until they move across to the outer part of the other wing. This outer part of the other wing in turn feeds back to the original wing. In this way, there is a natural dynamic disconnection into the red and blue sets seen in Figure \ref{LorenzRealDrift}.

\subsubsection{Optimally perturbing to decrease transport between the almost-invariant sets}

To decrease the transport between the almost-invariant sets, we perturb $P_0$ to maximally increase the magnitude of $\lambda_0$ so the almost-invariant sets exchange mass more slowly.
The rate of mixing between the two almost-invariant sets is controlled by the second eigenvalue $\lambda_0$. The tolerance is set to $\tau=10^{-12}$ as in the previous experiments. 
The optimal-response vector field is shown as arrows in Figure \ref{LorenzRealDrift}. 
All of the perturbation effort is spent pushing trajectories in the inner bands closer to the centres of rotation. 
This is because moving trajectories toward the rotation centres will mean those trajectories will remain longer in the corresponding almost-invariant set because it will take more time to spiral outwards.
No effort is put into the outer bands because those trajectories will shortly reenter the inner bands without additional perturbation.

\section{Case study 3: the El Ni\~no Southern Oscillation}
\label{sec:enso}

The El Ni\~no Southern Oscillation is the dominant mode of natural variability of the Earth's climate system on interannual (approximately 2--7 year) timescales.
It manifests itself as an SST oscillation in the eastern equatorial Pacific with an approximate 3--5 year periodicity, whose corresponding phases of anomalously warm and cold SSTs are known as El Ni\~no and La Ni\~na events, respectively.   
This oscillation modulates atmospheric circulation and weather patterns across the globe, generating considerable environmental and socioeconomic impacts. 

ENSO is an outcome of coupled nonlinear dynamics and thermodynamics involving the upper ocean and the atmosphere, the details of which have been a subject of active investigation over multiple decades \cite{TimmermannEtAl18}.   
One of the central research questions in this area is how ENSO is likely to respond to perturbations of the climate system, notably natural or anthropogenic changes in radiative forcing, e.g., due to volcanoes or greenhouse gas emissions.  
Optimal linear response provides a well-suited framework to characterize changes in the ENSO lifecycle under different perturbation scenarios using data from free-running equilibrium simulations of the climate. 

A standard approach for quantifying the strength and phase of the ENSO cycle is the Ni\~no~3.4 index, which is based on spatially-averaged SST anomalies over the box with coordinates [5$^\circ$N--5$^\circ$S, 170$^\circ$W--120$^\circ$W].
The Ni\~no~3.4 index is perhaps the most common index for monitoring ENSO due to the strong correlation between SST anomalies in the Ni\~no~3.4 box and ENSO activity \cite{BarnstonEtAl97}. 
It is typically computed as a moving mean over a window of 5 months, with El Ni\~no/La Ni\~na events declared when the index exceeds a threshold of $\pm 0.4\;\text{$^\circ$C}$ for a period of 6 consecutive months.   
In this work we use monthly-averaged Ni\~no 3.4 SST anomalies, without further moving averaging, to estimate the ENSO cycle and calculate optimal linear responses. We then illustrate these responses in terms of the SST anomalies and also in terms of  monthly-averaged Indo-Pacific SST and surface-wind anomaly fields, by considering them as observation maps in the optimal linear response vector field calculations.

In this section, we apply the techniques described in sections~\ref{sec:operator_approximation} and~\ref{sec:optlinresp} to identify optimal perturbations for modifying two key ENSO characteristics, namely mean period and decorrelation rate, using simulated data from a comprehensive climate model. 

\subsection{Extracting the ENSO cycle from sea-surface temperature anomaly data}

Our ENSO extraction approach closely follows \cite{FroylandEtAl21}.
The data is a time-ordered set of SST anomaly fields in an Indo-Pacific domain, taken from a 1300-year control integration of the Community Climate System Model Version 4 (CCSM4) \cite{GentEtAl11} with fixed pre-industrial forcings. Each SST field covers the region $[30^\circ \text{E}, 70^\circ \text{W}] \times [60^\circ \text{S}, 20^\circ \text{N}]$ with a mean longitudinal resolution of $0.36^\circ$ and mean latitudinal resolution of $0.41^\circ$. This region is sampled at $d = 44771$ points from the model's native ocean grid. The data is monthly averaged and sampled at a temporal frequency of one snapshot per month ($\Delta t = 1$ month) for $N=4800$ months covering the last 400 years of the 1300-year dataset. The seasonal cycle over this time span was removed from the raw SST fields, leaving only the anomalous SST components including the ENSO component.

To extract the ENSO cycle cleanly, the data is embedded via Takens delay-embedding with $Q = 12$ delays of length $L = 1$ month, resulting in an embedded time series of dimension $\tilde{d} = Qd = 537252$ and length $\tilde{N} = N - (Q-1)L = 4789$ months.
The stochastic transition matrix $P_0$ is constructed from the embedded trajectory as in \eqref{eq:markov_mat}, using $\epsilon = 0.03$, and with a step size of $s=6$, rather than the $s=1$ used in our previous case studies. Setting $s=6$ means replacing the indices $j-1$ and $k-1$ in \eqref{eq:markov_mat} with $j-6$ and $k-6$, respectively, to represent an evolution of $6\,\Delta t = 6$ months. Because $\epsilon$-noise is added at each step of the Markov chain, evolving in 6-month steps reduces the impact of short term (month-to-month) noise in the SST data. The eigenvalue $\lambda_0$ of $P_0$ also corresponds to a 6-month evolution, so to return to time units of months, from now on we replace all eigenvalues $\lambda$ with $\lambda^{1/6}$.

\begin{figure}
    \centering
\includegraphics[width=0.45\linewidth]{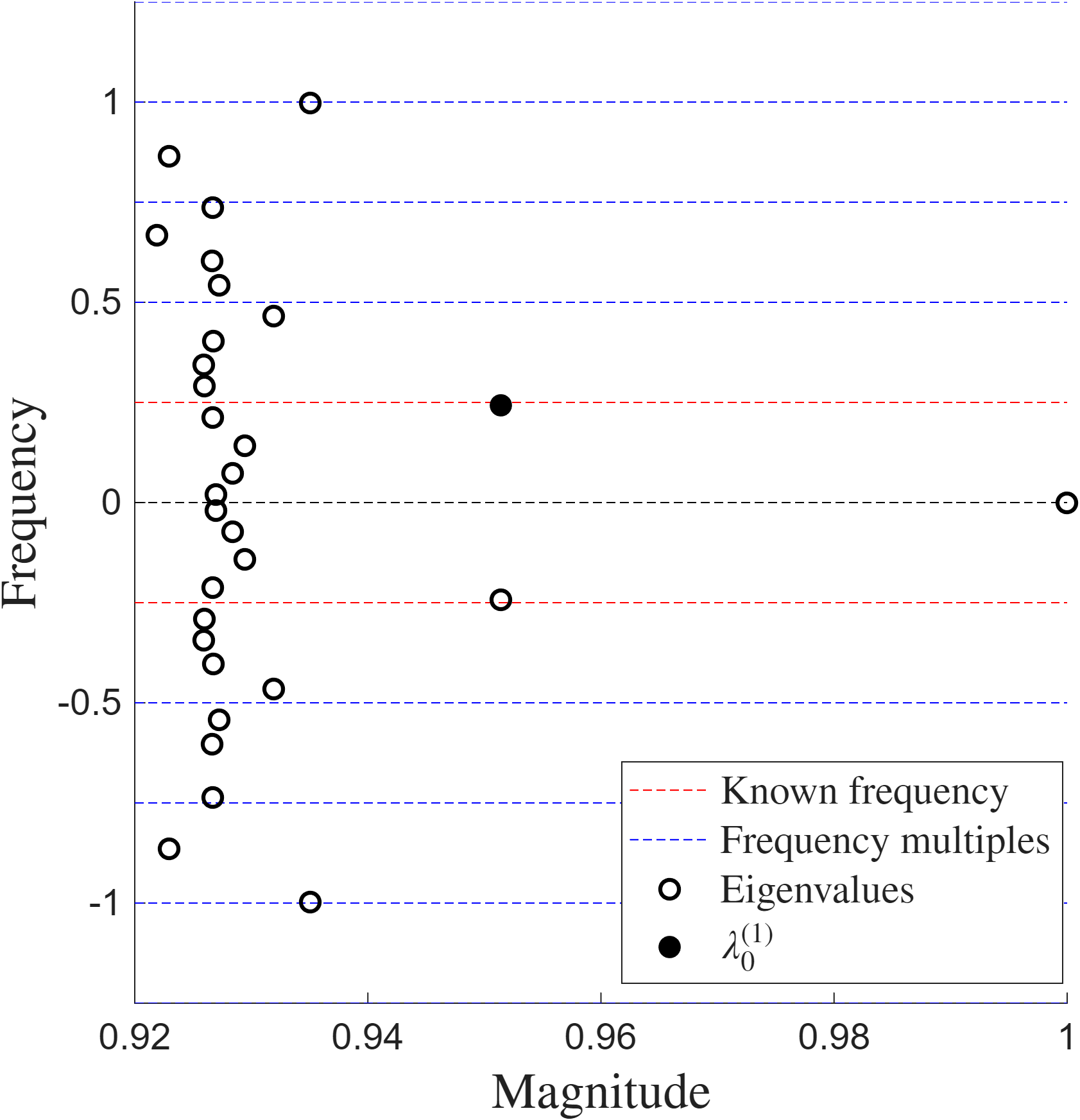}
\includegraphics[width=0.54\linewidth]{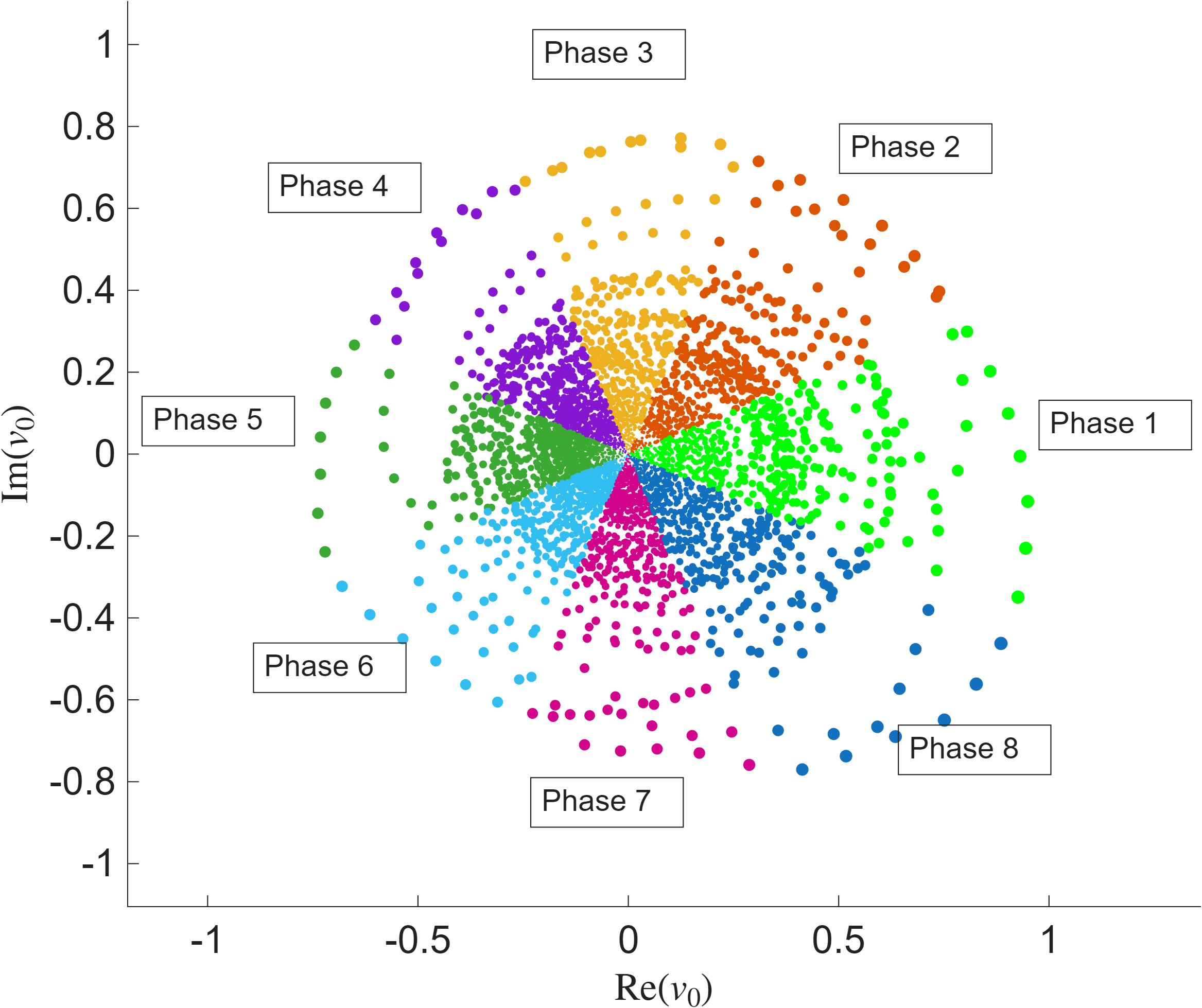}
    \caption{\emph{Left:} The magnitudes and frequencies of the largest 30 eigenvalues of $P_0$, approximating those of the transfer operator for the SST data. \emph{Right:} The real and imaginary parts of the ENSO eigenvector. The magnitude and angle of each point can be interpreted respectively as the ENSO strength and phase during each month of the data. The cycle is broken up into 8 phases (polar sectors) each spanning $\pi/4$ radians in the complex plane which corresponds to approximately 6 months.}
    \label{ENSOeig}
\end{figure}

Figure \ref{ENSOeig}(left) shows the spectrum of the unperturbed Markov matrix $P_0$. 
The leading nontrivial eigenvalues are $\lambda_0 = 0.9514 e^{i 0.1274}$ and its complex conjugate, which have a frequency of $\nu = \frac{\arg(\lambda_0)}{2\pi \Delta t} = 0.2433$ cycles per year or equivalently a period of 49.33 months. The large gap in magnitude to the remaining eigenvalues shows that ENSO is the strongest approximate cycle present in the dynamics. This gap is present for many choices of parameters $\epsilon, s, Q, L$, but seems to appear more consistently when $s\ge 6$ than with smaller integer step sizes. 
The values of the other parameters were chosen to increase the size of this gap and to ensure the corresponding eigenvector forms a roughly circular shape in $\mathbb{C}$ in Figure \ref{ENSOeig}(right).

When a complex eigenvector is computed with an eigensolver, it is typically $\ell^2$-normalised but has an arbitrary phase angle. Let $w_0$ be the eigenvector corresponding to $\lambda_0$ immediately after computation. We eliminate the arbitrary angle by setting $v_0 := w_0\exp(-i \arg((w_0)_k)$ where the $k$th month of the data is the one with the largest Ni\~no~3.4 signal, thus aligning zero phase with El Ni\~no. After this rotation, $v_0$ becomes a proxy indicator for ENSO phase and strength, with El Ni\~no events occurring when its real part is strongly positive and La Ni\~na occurring when strongly negative. Figure \ref{ENSOeig}(right) partitions $v_0$ into eight wedges, each occupied by the cycle for approximately 6 months based on the $\sim 49$ month periodicity inferred from the corresponding eigenvalue $\lambda_0$. The wedges traversed in the anticlockwise direction each represent a phase of the ENSO cycle, with Phase 1 being El Ni\~no and Phase 5 being La Ni\~na.
Details on the construction of the wedge sets and phase composites are given in section \ref{phasecompconstruct}.

\subsection{Construction of phase composites for interpretation of the ENSO eigenvector in observation space} \label{phasecompconstruct}

We recap the phase composite procedure from \cite{FroylandEtAl21}, which allows one to physically interpret the ENSO lifecycle captured by eigenvector $v_0$ and its associated optimal-response vector fields.
{A number $\mathcal{I} \in \mathbb{N}$ of phase composites are constructed} from equiangular sectors of the complex plane as follows.
Let $W_i\subset \{1,\ldots,\tilde{N}\}$, $i=1,\ldots,\mathcal{I}$, denote the set of indices of $v_0$ for which $(v_{0})_j$ lies in the $i^{\rm th}$ sector $\{\theta \in [\frac{(i-1-1/2)2\pi}{\mathcal{I}},\frac{(i-1+1/2)2\pi}{\mathcal{I}}) \mod 2\pi\}$.
Let $\bar{x}_j= x_{j+\lfloor (Q-1)L/2\rfloor}$ be the time-central element of the delay-embedded data point $\tilde{x}_j=[x_j,x_{j+L},\dots,x_{j+(Q-1)L}] \in \mathbb{R}^{Qd}$.

For each $i=1,\dots,\mathcal{I}$ we construct phase composites as eigenvector-magnitude-weighted convex combinations of time-central data points restricted to the $i^{\rm th}$ sector.
\begin{equation} \label{phasecomps}
    \chi_i \coloneq \frac{\sum_{j\in{W_i}}|(v_0)_j|\bar{x}_j}{\sum_{j\in{W_i}}|(v_0)_j|}\in\mathbb{R}^d.
\end{equation}
When the eigenvector has a large magnitude, the corresponding observations will be weighted more heavily because they are more strongly associated with the cycle described by $v_0$ and $\lambda_0$.
In \cite{FroylandEtAl21} similar phase composites were constructed with \eqref{phasecomps}, replacing the weighting $x\mapsto |x|$ with $x\mapsto \mathbf{1}_{|x|>r}(x)$ for some cut-off radius $0<r<1$.

 In Figure \ref{ENSOeig}(right) we partition the eigenvector into $\mathcal{I} = 8$ sectors, each representing a $49.33/8 = 6.169$ month phase of the cycle. {Taking $\bar{x}_j$ as the time-central element of the time-delayed SST field $\tilde{x}_j=F_{Q,L}(\omega_j)$, we construct the 8-phase composite SST fields $\chi_i$, shown in the first column of Figure \ref{ENSOPhaseMaps}. Each $\chi_i$ is an SST anomaly field representative of typical observations from Phase $i$ of the ENSO cycle.}

As noted in section \ref{sec:phase_space_perturbations}, we are also able to construct phase composites $\chi'_i \in \mathbb{R}^{d'}$ for a different set of observations given by a function $F':M\rightarrow \mathbb{R}^{d'}$, by replacing the central element $\bar{x}_j$ of $\tilde{x}_j=F_{Q,L}(\omega_j)$ in \eqref{phasecomps} with $\bar{x}'_j$, the central element of $\tilde{x}_j'=F_{Q,L}'(\omega_j)$. {To aid the physical interpretation of the phases, we also construct phase composites for the wind-speed anomaly fields, by replacing the SST data $\bar{x}_j$ (observed with $F$) in \eqref{phasecomps} with surface wind data $\bar{x}_j'$ (observed with $F'$).} {This next sentence I cut and pasted from section 6.3:}
The surface wind data is collected from the same CCSM4 simulation as the SST data, and consists of {zonal (West--East) and meridional (North--South) atmospheric surface wind components} at each point of a grid with mean longitudinal resolution of $1.25^\circ$ and a mean latitudinal resolution of $0.94^\circ$.

\subsection{Phase composites of the optimal-response vector field}

To construct phase composites for the optimal-response vector fields we create the $\tilde{N}\times d$ array $\bar{X}$ whose $j^{\rm th}$ row is  $\bar{x}_j\in\mathbb{R}^d$, and modify \eqref{phasecomps} as:
\begin{equation} \label{phaseresponse}
    \dot{\chi}_i = \frac{\sum_{j\in{W_i}}|(v_0)_j|(\dot{P}\bar{X})_j}{\sum_{j\in{W_i}}|(v_0)_j|}\in\mathbb{R}^d.
\end{equation}
Each $\dot{\chi}_i$, $i=1,\ldots,\mathcal{I}$, represents the  response of the observable $F$ to the perturbation $\dot{P}$ in the $i$th sector of the cycle.
Replacing $\bar{X}$ with $\bar{X}'$ in \eqref{phaseresponse}, similar to the replacement discussed at the end of the previous subsection, yields the phase composites for the optimal-response vector field associated with the observable $F'$ instead.

\subsection{Perturbations of the ENSO cycle}
\label{ENSOPerturb}

We compute perturbation matrices $\dot{P}$, to (i) maximally increase the frequency of the ENSO cycle and to (ii) maximally increase the temporal correlation strength of the ENSO cycle. 
For both perturbations, we construct $\dot{P}$ with a tolerance of $\tau=3 \times 10^{-4}$, resulting in $\dot{P}$ being $9\%$ nonzero.

\subsubsection{Optimally decreasing the period of the ENSO cycle}
\label{sec:ensoperiod}

\begin{figure}
    \centering
\includegraphics[width=\linewidth]{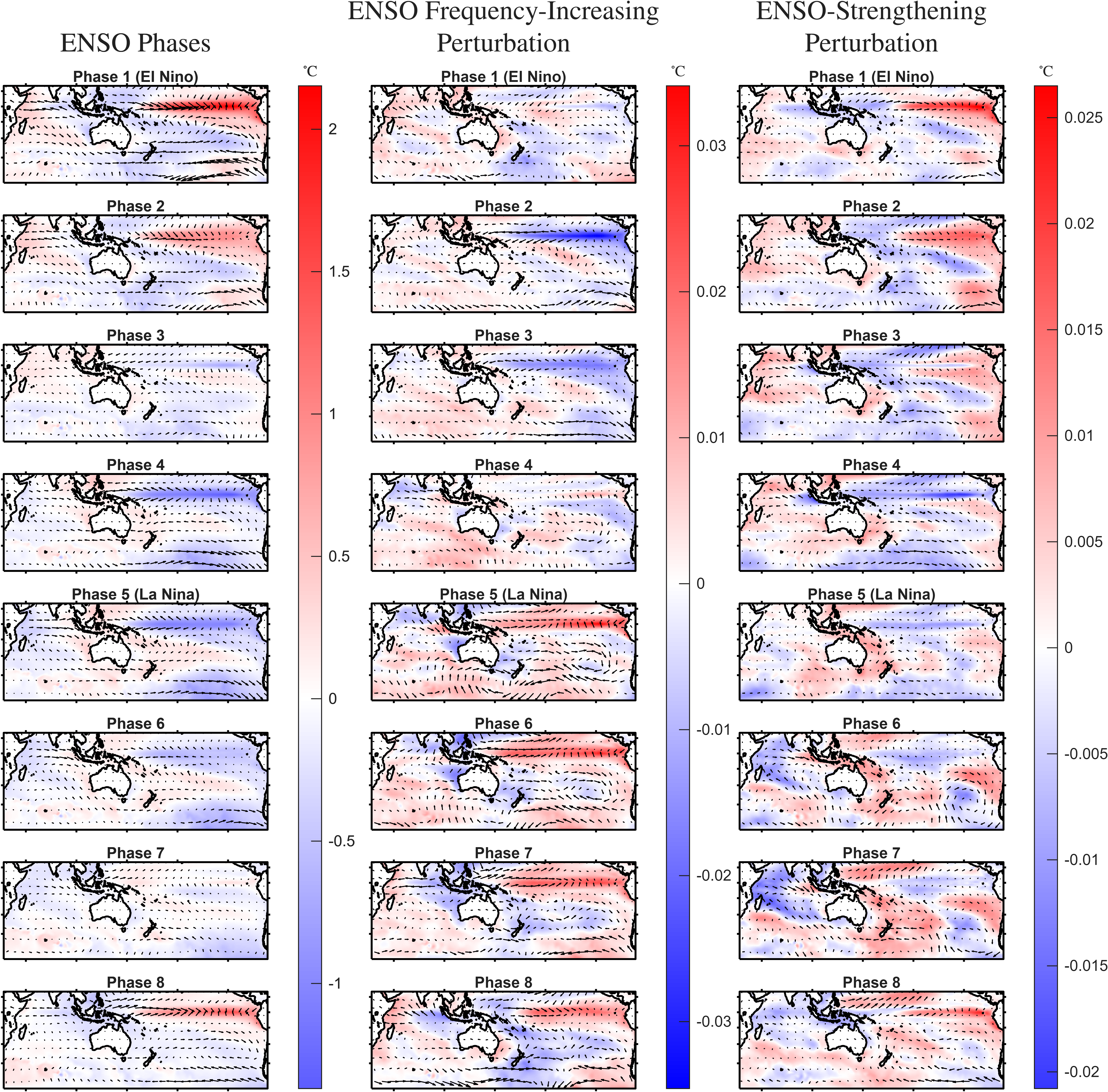}
    \caption{\emph{Left:} Phase composites of SST anomalies (colours) and surface wind anomalies (arrows) based on the ENSO eigenvector $v_0$. Each arrow represents the mean wind anomaly in the surrounding $7.5^\circ \times 7.5^\circ$ square. \emph{Center:} Phase composites of the optimal-response vector fields for the SST field (colours) and surface wind field (arrows), for the ENSO frequency-increasing perturbation. \emph{Right:} The same information for the ENSO-strengthening perturbation.}
    \label{ENSOPhaseMaps}
\end{figure}

For the frequency-increasing optimisation problem, the optimal value of the derivative of $\arg \lambda_\delta$ in \eqref{eq:freqderivative} is 0.5132. To physically interpret our optimal frequency perturbations, we compute the optimal-response vector fields described in section \ref{sec:phase_space_perturbations} for both the SST and surface wind data sets, and then use \eqref{phaseresponse} to construct the phase composites of the response fields. For both observables, the response phase composites for the frequency-increasing perturbation are displayed in the second column of Figure~\ref{ENSOPhaseMaps}.
We see in the second column of Figure \ref{ENSOPhaseMaps} that the frequency perturbation is \emph{leading} the ENSO cycle shown in the first column of Figure \ref{ENSOPhaseMaps}, with an approximate phase difference of $\pi/2$.
Adding SST and anomalous surface wind with a $\pi/2$ phase difference (i.e., during the neutral phases of the cycle) will tend to shorten the duration of El Ni\~no to La Ni\~na transitions, thus increasing the frequency of the cycle.
Away from the Ni\~no 3.4 regions, it is worthwhile noting that the response vector fields in Figure~\ref{ENSOPhaseMaps} exhibit significant activity in the tropical Indian Ocean (e.g., during Phase~7).

\subsubsection{Optimally strengthening the ENSO cycle}
\label{sec:ensostrength}

For the eigenvalue magnitude-increasing optimisation problem, the optimal value of the derivative \eqref{eq:magderiv} is 0.5098, which corresponds to a derivative in the magnitude of $\lambda_0$ of 0.3791. The optimal-response vector field to increase the temporal correlation strength of the ENSO cycle is seen in the third column of Figure \ref{ENSOPhaseMaps} to be \emph{in phase} with the underlying ENSO cycle.
The in-phase nature of the magnitude perturbation, i.e. adding SST and anomalous surface wind sources coherently with the ENSO cycle will also tend to increase its amplitude.

The IOD \cite{HameedEtAl99,SajiEtAl99,WebsterEtAl99} is a prominent pattern of interannual variability in the Indian Ocean that is known to exhibit non-trivial correlations with ENSO \cite{SlawinskaGiannakis17, HameedEtAl18}. 
Its positive phases are characterised by warmer SSTs in the western Indian Ocean and cooler in the east, which are reversed during its negative phase \cite{SajiEtAl99,WebsterEtAl99}. 
A positive IOD pattern is present in the third column of Figure~\ref{ENSOPhaseMaps} in phases 3 and 4, and a negative IOD pattern is present in phases 6 and 7.
This suggests that enhancing the predictability of the ENSO cycle involves some interaction with the IOD: more specifically, a negative IOD state during El Ni\~no buildup and a positive IOD state during La Ni\~na buildup. 
Previous work has found that a positive IOD signal preceding an El Ni\~no increases the likelihood of so-called super El Ni\~no events \cite{HameedEtAl18}.
Our results in the third column of Figure \ref{fig:nino_34_composites} concern ENSO predictability-enhancing perturbations (as opposed to amplitude increasing mechanisms), and show a different IOD interaction.

\subsubsection{Spatially averaged optimal-response vector fields in the Ni\~no~3.4 box}

\begin{figure}
    \centering
\includegraphics[width=0.5\linewidth]{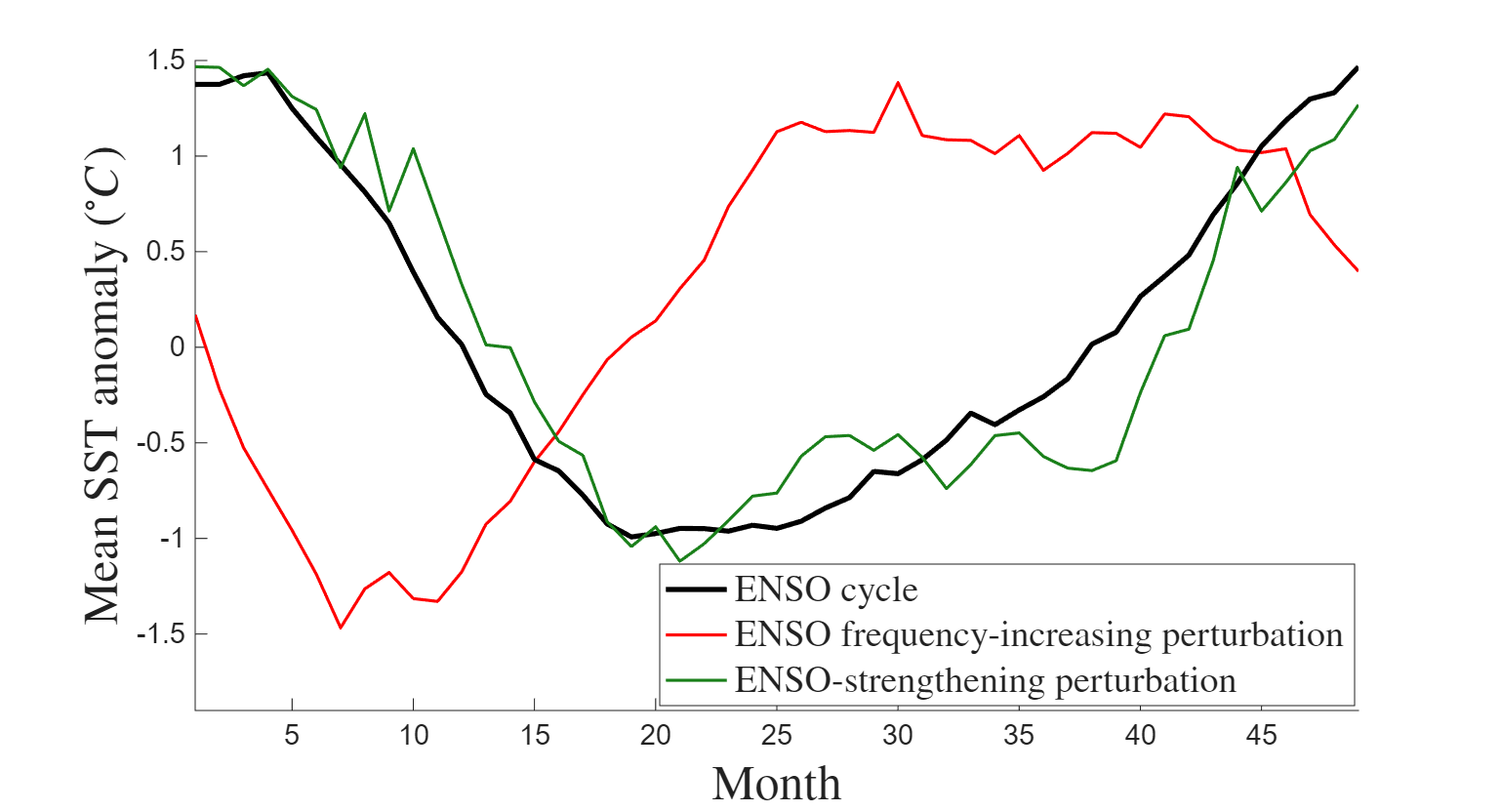}
    \caption{SST anomalies (black line) and optimal-response SST fields for the frequency-increasing perturbation (red line) and magnitude-increasing perturbation (green line), averaged over the Ni\~no 3.4 box for each month of the ENSO cycle.}
    \label{fig:nino_34_composites}
\end{figure}

According to the argument of $\lambda_0$, the ENSO cycle has a period of 49.33 months. For each month, Figure \ref{fig:nino_34_composites} shows the mean temperature anomaly in the Ni\~no~3.4 box, computed using the phase composite method \eqref{phasecomps}. To account for the non-integer period, we partition the eigenvector $v_0$ into $\lceil 49.33 \rceil = 50$ sectors of the complex plane, with the first 49 sectors being equiangular and the 50th covering a smaller angle that is $33\%$ the size of the other sectors. We construct 50 SST phase composites using  \eqref{phasecomps} and the sets $W_i\subset \{1,\ldots,\tilde{N}\}$, $i=1,\ldots,50$. Finally, to produce Figure~\ref{fig:nino_34_composites} we calculate the mean temperature anomaly within the Ni\~no~3.4 box for each of the 50 phase composite SST maps.

The computations in sections \ref{sec:ensoperiod} and \ref{sec:ensostrength} use the disallowed transition set $D_\tau = \{(i,j) \in \{1,\dots, \tilde{N}\}:\tau<P_{ij}<1-\tau\}$ with $\tau = 3 \times 10^{-4}$, resulting in perturbation matrices $\dot{P}$ that are $9\%$ nonzero.
This means that the optimisation can only adjust the conditional probabilities of transitions from $\tilde{x}_i$ to the future embedded states $\tilde{x}_j$ that are the 9\%-closest to $\tilde{x}_{i+1}$, thus preferring ``small, physical perturbations'' in a  small neighbourhood of SST one-year histories. Enlarging this neighbourhood by decreasing $\tau$ allows for more degrees of freedom in maximising $\left.\frac{\partial |\lambda_\delta|}{\partial \delta}\right|_{\delta = 0}$ and $\left.\frac{\partial \arg\lambda_\delta}{\partial \delta}\right|_{\delta = 0}$, at the expense of using larger and potentially less physical perturbations. With this background, we can report that the results in Figures \ref{ENSOPhaseMaps} and \ref{fig:nino_34_composites} remain relatively unchanged when we reduce $\tau$ to $10^{-4}$, yielding perturbation matrices that are $26\%$ nonzero. 
Thus, our optimal perturbations are robust to moderate changes in $\tau$.

\section{Conclusions}
\label{sec:conclusions}
We developed data-driven techniques for computation of optimal linear response, extending the results of \cite{AntownFroyland18} to {handle high-dimensional time-series data arising from observations.}
This is achieved using kernel smoothing of the transfer/Koopman operator of the dynamical system to produce a finite-state Markov chain approximation. The optimisation problem for manipulating the dominant spectrum (corresponding to almost-cycles and almost-invariant sets) of the Koopman operator is discretised to an optimal perturbation problem for the Markov chain approximation, making its solution computationally tractable.
We also introduce a notion of optimal-response vector fields, computed from the Markov chain perturbation, as a tool for visualising and physically interpreting the optimal response.

{In our numerical experiments, we focused on perturbations that maximally increase the frequency or suppress the decay of correlations of dominant almost-cycles, or suppress decay of almost-invariant sets.} {We illustrated the efficacy of the optimal-response vector field} across several case studies, chiefly the Lorenz system and the ENSO cycle.
The frequency-increasing optimal-response vector field for the Lorenz system displayed an intricate mechanism, pointing orthogonal to level sets of the argument of the complex eigenvector, instead of pointing in the direction of the flow.
Meanwhile, the ENSO case study revealed a connection between SST perturbations that enhance ENSO cycle predictability and a specific phasing of the IOD pattern.

More broadly, apart from the response of the transfer/Koopman operator spectrum studied in this paper, our numerical techniques may be directly applied to optimise general linear responses.
Future work will pursue these applications to answer a variety scientific questions.

\section*{Acknowledgements}

The research of GF is supported by an Australian Research Council Laureate Fellowship (FL230100088) and partially supported by an Einstein Visiting Fellowship (Einstein Foundation Berlin).
{DG acknowledges support from the U.S.\ Department of Defense, Basic Research Office, under Vannevar Bush Faculty Fellowship grant N00014-21-1-2946 and the U.S.\ Office of Naval Research under MURI grant N00014-19-1-242.}
{The research of NP was supported by the Commonwealth through an Australian Government Research Training Program Scholarship [DOI: https://doi.org/10.82133/C42F-K220], and partially supported through an ARC Laureate Fellowship.}


\begin{thebibliography}{89}
\providecommand{\natexlab}[1]{#1}
\providecommand{\url}[1]{\texttt{#1}}
\expandafter\ifx\csname urlstyle\endcsname\relax
  \providecommand{\doi}[1]{doi: #1}\else
  \providecommand{\doi}{doi: \begingroup \urlstyle{rm}\Url}\fi

\bibitem[Ghil and Lucarini(2020)]{GhilLucarini20}
M.~Ghil and V.~Lucarini.
\newblock The physics of climate variability and climate change.
\newblock \emph{Rev. Mod. Phys.}, 92:\penalty0 035002, 2020.
\newblock \doi{10.1103/RevModPhys.92.035002}.

\bibitem[Lucarini et~al.(2017)Lucarini, Ragone, and Lunkeit]{LucariniEtAl17}
V.~Lucarini, F.~Ragone, and F.~Lunkeit.
\newblock Predicting climate change using response theory: {G}lobal averages
  and spatial patterns.
\newblock \emph{J. Stat. Phys.}, 166:\penalty0 1036--1064, 2017.
\newblock \doi{10.1007/s10955-016-1506-z}.

\bibitem[Ruelle(1997)]{Ruelle97}
D.~Ruelle.
\newblock Differentiation of {SRB} states.
\newblock \emph{Commun. Math. Phys.}, 187:\penalty0 227--241, 1997.
\newblock \doi{10.1007/s002200050134}.

\bibitem[Ruelle(2009)]{Ruelle09}
D.~Ruelle.
\newblock A review of linear response theory for general differentiable
  dynamical systems.
\newblock \emph{Nonlinearity}, 22\penalty0 (4):\penalty0 855--870, 2009.
\newblock \doi{10.1088/0951-7715/22/4/009}.

\bibitem[Baladi(2014)]{Baladi14}
V.~Baladi.
\newblock Linear response, or else, 2014.
\newblock URL \url{https://arxiv.org/abs/1408.2937}.
\newblock ICM Seoul 2014 talk.

\bibitem[Callen and Welton(1951)]{CallenWelton51}
H.~B. Callen and T.~A. Welton.
\newblock Irreversibility and generalized noise.
\newblock \emph{Phys. Rev.}, 83\penalty0 (1):\penalty0 34--40, 1951.
\newblock \doi{10.1103/PhysRev.83.34}.

\bibitem[Kubo(1966)]{Kubo66}
R.~Kubo.
\newblock The fluctuation--dissipation theorem.
\newblock \emph{Rep. Prog. Phys.}, 29:\penalty0 255--284, 1966.
\newblock \doi{10.1088/0034-4885/29/1/306}.

\bibitem[Kraichnan(1959)]{Kraichnan59}
R.~H. Kraichnan.
\newblock Classical fluctuation--relaxation theorem.
\newblock \emph{Phys. Rev.}, 113\penalty0 (5), 1959.
\newblock \doi{10.1103/PhysRev.113.1181}.

\bibitem[Leith(1975)]{Leith75}
C.~E. Leith.
\newblock Climate response and fluctuation dissipation.
\newblock \emph{J. Atmos. Sci.}, 32\penalty0 (10):\penalty0 2022--2026, 1975.
\newblock \doi{https://doi.org/10.1175/1520-0469(1975)032<2022:CRAFD>2.0.CO;2}.

\bibitem[Gallavotti and Cohen(1995)]{GallavottiCohen95}
G.~Gallavotti and E.~G.~D. Cohen.
\newblock Dynamic ensembles in stationary states.
\newblock \emph{Phys. Rev. Lett.}, 74\penalty0 (14):\penalty0 2694--2697, 1995.
\newblock \doi{10.1103/PhysRevLett.74.2694}.

\bibitem[Gallavotti(1996)]{Gallavotti96}
G.~Gallavotti.
\newblock Onsager hypothesis: {O}nsager reciprocity and
  fluctuation--dissipation theorem.
\newblock \emph{J. Stat. Phys.}, 84\penalty0 (5/6):\penalty0 899--925, 1996.
\newblock \doi{10.1007/bf02174123}.

\bibitem[Young(2002)]{Young02}
L.-S. Young.
\newblock What are {SRB} measures, and which dynamical systems have them?
\newblock \emph{J. Stat. Phys.}, 108:\penalty0 733--754, 2002.
\newblock \doi{10.1023/a:1019762724717}.

\bibitem[Bowen and Ruelle(1975)]{BowenRuelle75}
R.~Bowen and D.~Ruelle.
\newblock The ergodic theory of {A}xiom {A} flows.
\newblock \emph{Inventiones Mathematicae}, 29:\penalty0 181--202, 1975.
\newblock \doi{10.1007/bf01389848}.

\bibitem[Gou{\"e}zel and Liverani(2006)]{GouezelLiverani06}
S.~Gou{\"e}zel and C.~Liverani.
\newblock Banach spaces adapted to {A}nosov systems.
\newblock \emph{Ergod. Theory Dyn. Syst.}, 26:\penalty0 189--217, 2006.
\newblock \doi{10.1017/s0143385705000374}.

\bibitem[Abramov and Majda(2008)]{AbramovMajda08}
R.~V. Abramov and A.~Majda.
\newblock New approximations and tests of linear fluctuation-response for
  chaotic nonlinear forced-dissipative dynamical systems.
\newblock \emph{J. Nonlinear Sci.}, 18:\penalty0 303--341, 2008.
\newblock \doi{10.1007/s00332-007-9011-9}.

\bibitem[Baladi and Smania(2009)]{BaladiSmania09}
V.~Baladi and D.~Smania.
\newblock Smooth deformations of piecewise expanding unimodal maps.
\newblock \emph{Discrete and Continuous Dynamical Systems}, 23\penalty0
  (3):\penalty0 685–703, 2009.
\newblock \doi{10.3934/dcds.2009.23.685}.

\bibitem[Hairer and Majda(2010)]{HairerMajda10}
M.~Hairer and A.~J. Majda.
\newblock A simple framework to justify linear response theory.
\newblock \emph{Nonlinearity}, 23\penalty0 (4):\penalty0 909--922, 2010.
\newblock \doi{10.1088/0951-7715/23/4/008}.

\bibitem[Baladi et~al.(2017)Baladi, Kuna, and Lucarini]{BaladiKunaLucarini17}
V.~Baladi, T.~Kuna, and V.~Lucarini.
\newblock Linear and fractional response for the srb measure of smooth
  hyperbolic attractors and discontinuous observables.
\newblock \emph{Nonlinearity}, 30\penalty0 (3), 2017.
\newblock \doi{10.1088/1361-6544/aa5b13}.

\bibitem[Baladi and Todd(2016)]{BaladiTodd16}
V.~Baladi and M.~Todd.
\newblock Linear response for intermittent maps.
\newblock \emph{Communications in Mathematical Physics}, 347, 2016.
\newblock \doi{10.1007/s00220-016-2577-z}.

\bibitem[Wormell and Gottwald(2018)]{WormellGottwald18}
C.L. Wormell and G.A. Gottwald.
\newblock On the validity of linear response theory in high-dimensional
  deterministic dynamical systems.
\newblock \emph{Journal of Statistical Physics}, 172, 2018.
\newblock \doi{10.1007/s10955-018-2106-x}.

\bibitem[{Gottwald} et~al.(2016){Gottwald}, Wormell, and
  Wouters]{GottwaldWormellWouters16}
G.~A. {Gottwald}, C.~Wormell, and J~Wouters.
\newblock On spurious detection of linear response and misuse of the
  fluctuation-dissipation theorem in finite time series.
\newblock \emph{Physica D Nonlinear Phenomena}, 331:\penalty0 89--101, 2016.
\newblock \doi{10.1016/j.physd.2016.05.010}.

\bibitem[Bahsoun et~al.(2018)Bahsoun, Galatolo, Nisoli, and
  Niu]{bahsoun2018rigorous}
Wael Bahsoun, Stefano Galatolo, Isaia Nisoli, and Xiaolong Niu.
\newblock A rigorous computational approach to linear response.
\newblock \emph{Nonlinearity}, 31\penalty0 (3):\penalty0 1073--1109, 2018.

\bibitem[Chandramoorthy and Wang(2022)]{ChandramoorthyWang22}
N.~Chandramoorthy and Q.~Wang.
\newblock Efficient computation of linear response of chaotic attractors with
  one-dimensional unstable manifolds.
\newblock \emph{SIAM J. Appl. Dyn. Syst.}, 21\penalty0 (2):\penalty0 735--781,
  2022.
\newblock \doi{10.1137/21m1405599}.

\bibitem[Ni(2026)]{ni2026fast}
A.~Ni.
\newblock Fast differentiation of hyperbolic chaos.
\newblock \emph{Archive for Rational Mechanics and Analysis}, 250\penalty0 (1),
  2026.

\bibitem[Froyland and Phalempin(2025)]{FroylandPhalempin25}
G.~Froyland and M.~Phalempin.
\newblock Optimal linear response for {A}nosov diffeomorphisms, 2025.
\newblock URL \url{https://arxiv.org/abs/2504.16532}.

\bibitem[Bell(1980)]{Bell80}
T.~Bell.
\newblock Climate sensitivity from fluctuation dissipation: {S}ome simple model
  tests.
\newblock \emph{J. Atmos. Sci}, 37\penalty0 (8):\penalty0 1700--1708, 1980.
\newblock \doi{10.1175/1520-0469(1980)037<1700:csffds>2.0.co;2}.

\bibitem[Gritsoun(2001)]{Gritsoun01}
A.~S. Gritsoun.
\newblock Fluctuation-dissipation theorem on attractors of atmospheric models.
\newblock \emph{Russ. J. Numer. Anal. Math. Modelling}, 16\penalty0
  (2):\penalty0 115--133, 2001.
\newblock \doi{10.1515/rnam-2001-0203}.

\bibitem[Gritsoun et~al.(2002)Gritsoun, Branstator, and
  Dymnikov]{GritsounEtAl02}
A.~S. Gritsoun, G.~Branstator, and V.~P. Dymnikov.
\newblock Construction of the linear response operator of an atmospheric
  general circulation model to small external forcing.
\newblock \emph{Russ. J. Numer. Anal. Math. Modelling}, 17\penalty0
  (5):\penalty0 399--416, 2002.
\newblock \doi{10.1515/rnam-2002-0503}.

\bibitem[Majda et~al.(2005)Majda, Abramov, and Grote]{MajdaEtAl05}
A.~J. Majda, R.~V. Abramov, and M.~J. Grote.
\newblock \emph{Information Theory and Stochastics for Multiscale Nonlinear
  Systems}, volume~25 of \emph{{CRM} Monograph Series}.
\newblock Americal Mathematical Society, Providence, 2005.

\bibitem[Majda et~al.(2010)Majda, Abramov, and Gershgorin]{MajdaEtAl10b}
A.~J. Majda, R.~Abramov, and B.~Gershgorin.
\newblock High skill in low-frequency climate response through fluctuation
  dissipation theorems despite structural instability.
\newblock \emph{Proc. Natl. Acad. Sci}, 107\penalty0 (2):\penalty0 581--586,
  2010.
\newblock \doi{10.1073/pnas.091299710}.

\bibitem[Dellnitz and Junge(1999)]{DellnitzJunge99}
M.~Dellnitz and O.~Junge.
\newblock On the approximation of complicated dynamical behavior.
\newblock \emph{SIAM J. Numer. Anal.}, 36:\penalty0 491, 1999.
\newblock \doi{10.1137/s0036142996313002}.

\bibitem[Sch\"utte et~al.(2001)Sch\"utte, Huisinga, and
  Deuflhard]{SchutteEtAl01}
Ch. Sch\"utte, W.~Huisinga, and P.~Deuflhard.
\newblock Transfer operator approach to conformational dynamics in biomolecular
  systems.
\newblock In B.~Fiedler, editor, \emph{Ergodic Theory, Analysis, and Efficient
  Simulation of Dynamical Systems}, pages 191--223. Springer-Verlag, Berlin,
  2001.
\newblock \doi{10.1007/978-3-642-56589-2_9}.

\bibitem[Mezi\'c(2005)]{Mezic05}
I.~Mezi\'c.
\newblock Spectral properties of dynamical systems, model reduction and
  decompositions.
\newblock \emph{Nonlinear Dyn.}, 41:\penalty0 309--325, 2005.
\newblock \doi{10.1007/s11071-005-2824-x}.

\bibitem[Mezi\'c and Banaszuk(2004)]{MezicBanaszuk04}
I.~Mezi\'c and A.~Banaszuk.
\newblock Comparison of systems with complex behavior.
\newblock \emph{Phys. D.}, 197:\penalty0 101--133, 2004.
\newblock \doi{10.1016/j.physd.2004.06.015}.

\bibitem[Otto and Rowley(2021)]{OttoRowley21}
S.~E. Otto and C.~W. Rowley.
\newblock Koopman operators for estimation and control of dynamical systems.
\newblock \emph{Annu. Rev. Control Robot. Auton. Syst.}, 4:\penalty0 59--87,
  2021.
\newblock \doi{10.1146/annurev-control-071020-010108}.

\bibitem[Colbrook(2024)]{Colbrook24}
M.~Colbrook.
\newblock The multiverse of dynamic mode decomposition algorithms.
\newblock In \emph{Handbook of Numerical Analysis}, pages 127--230. Amsterdam,
  2024.
\newblock \doi{10.1016/bs.hna.2024.05.004}.

\bibitem[Brunton et~al.(2022)Brunton, Budisi{\'c}, Kaiser, and
  Kutz]{BruntonEtAl22}
S.~L. Brunton, M.~Budisi{\'c}, E.~Kaiser, and J.~N. Kutz.
\newblock Modern {K}oopman theory for dynamical systems.
\newblock \emph{SIAM Rev.}, 64\penalty0 (2):\penalty0 229--340, 2022.
\newblock \doi{10.1137/21m1401243}.

\bibitem[Baladi(2018)]{Baladi2018}
V.~Baladi.
\newblock \emph{Dynamical Zeta Functions and Dynamical Determinants for
  Hyperbolic Maps: A Functional Approach}.
\newblock Springer, 2018.
\newblock \doi{https://doi.org/10.1007/978-3-319-77661-3}.

\bibitem[Santos~Guti{\'e}rrez and Lucarini(2022)]{SantosGutierrezLucarini22}
M.~Santos~Guti{\'e}rrez and V.~Lucarini.
\newblock On some aspects of the response to stochastic and deterministic
  forcings.
\newblock \emph{J. Phys. A: Math. Theor.}, art. 425002, 2022.
\newblock \doi{10.1088/1751-8121/ac90fd}.

\bibitem[Zagli et~al.(2026)Zagli, Colbrook, Lucarini, Mezi{'c}, and
  Moroney]{ZagliEtAl26}
N.~Zagli, M.~J. Colbrook, V.~Lucarini, I.~Mezi{'c}, and J.~Moroney.
\newblock Bridging the gap between {K}oopmanism and response theory: {U}sing
  natural variability to predict forced response.
\newblock \emph{SIAM J. Appl. Dyn. Syst.}, 25\penalty0 (1):\penalty0 196--229,
  2026.
\newblock \doi{10.1137/24m1699206}.

\bibitem[Lucarini et~al.(2026)Lucarini, Santos~Guti{\'e}rrez, Moroney, and
  Zagli]{LucariniEtAl26}
V.~Lucarini, M.~Santos~Guti{\'e}rrez, J.~Moroney, and N.~Zagli.
\newblock A general framework for linking free and forced fluctuations via
  {K}oopmanism.
\newblock \emph{Chaos Solitons Fractals}, 202\penalty0 (Part 1):\penalty0
  117540, 2026.
\newblock \doi{10.1016/j.chaos.2025.117540}.

\bibitem[Williams et~al.(2015)Williams, Kevrekidis, and Rowley]{WilliamsEtAl15}
M.~O. Williams, I.~G. Kevrekidis, and C.~W. Rowley.
\newblock A data-driven approximation of the {K}oopman operator: Extending
  dynamic mode decomposition.
\newblock \emph{J. Nonlinear Sci.}, 25\penalty0 (6):\penalty0 1307--1346, 2015.
\newblock \doi{10.1007/s00332-015-9258-5}.

\bibitem[Dellnitz et~al.(2000)Dellnitz, Froyland, and Sertl]{DellnitzEtAl00}
M.~Dellnitz, G.~Froyland, and S.~Sertl.
\newblock On the isolated spectrum of the {P}erron–{F}robenius operator.
\newblock \emph{Nonlinearity}, 13:\penalty0 1171--1188, 2000.
\newblock \doi{10.1088/0951-7715/13/4/310}.

\bibitem[Froyland et~al.(2014)Froyland, Gonz\'alez-Tokman, and
  Quas]{FroylandEtAl14b}
G.~Froyland, C.~Gonz\'alez-Tokman, and A.~Quas.
\newblock Detecting isolated spectrum of transfer and {K}oopman operators with
  {F}ourier analytic tools.
\newblock \emph{J. Comput. Dyn.}, 1\penalty0 (2):\penalty0 249--278, 2014.
\newblock \doi{10.3934/jcd.2014.1.249}.

\bibitem[Porte(2019)]{Porte19}
M.~Porte.
\newblock Linear response for {D}irac observables of {A}nosov diffeomorphisms.
\newblock \emph{Discrete and Continuous Dynamical Systems}, 39\penalty0
  (4):\penalty0 1799--1819, 2019.
\newblock \doi{10.3934/dcds.2019078}.

\bibitem[Froyland(2005)]{Froyland05}
G.~Froyland.
\newblock Statistically optimal almost-invariant sets.
\newblock \emph{Phys. D.}, 200:\penalty0 205--219, 2005.
\newblock \doi{10.1016/j.physd.2004.11.008}.

\bibitem[Froyland et~al.(2021)Froyland, Giannakis, Lintner, Pike, and
  Slawinska]{FroylandEtAl21}
G.~Froyland, D.~Giannakis, B.~Lintner, M.~Pike, and J.~Slawinska.
\newblock Spectral analysis of climate dynamics with operator-theoretic
  approaches.
\newblock \emph{Nat. Commun.}, 12:\penalty0 6570, 2021.
\newblock \doi{10.1038/s41467-021-26357-x}.

\bibitem[Castro and Froyland(2025)]{CastroFroyland25}
M.~M. Castro and G.~Froyland.
\newblock On the structure of complex spectra and eigenfunctions of transfer
  and {K}oopman operators, 2025.
\newblock URL \url{https://arxiv.org/abs/2505.05770}.

\bibitem[Froyland and Santitissadeekorn(2017)]{froyland2017optimal}
G.~Froyland and N.~Santitissadeekorn.
\newblock Optimal mixing enhancement.
\newblock \emph{SIAM Journal on Applied Mathematics}, pages 1444--1470, 2017.

\bibitem[Froyland et~al.(2020)Froyland, Koltai, and
  Stahn]{froyland2020computation}
G.~Froyland, P.~Koltai, and M.~Stahn.
\newblock Computation and optimal perturbation of finite-time coherent sets for
  aperiodic flows without trajectory integration.
\newblock \emph{SIAM Journal on Applied Dynamical Systems}, 19\penalty0
  (3):\penalty0 1659--1700, 2020.

\bibitem[Antown et~al.(2018)Antown, Dragi\v{c}evi\'{c}, and
  Froyland]{AntownFroyland18}
F.~Antown, D.~Dragi\v{c}evi\'{c}, and G.~Froyland.
\newblock Optimal linear responses for {M}arkov chains and stochastically
  perturbed dynamical systems.
\newblock \emph{J. Stat. Phys.}, 170\penalty0 (6):\penalty0 1051–1087, 2018.
\newblock \doi{10.1007/s10955-018-1985-1}.

\bibitem[Antown et~al.(2022)Antown, Froyland, and Galatolo]{AntownEtAl22}
F.~Antown, G.~Froyland, and S.~Galatolo.
\newblock Optimal linear response for {M}arkov {H}ilbert--{S}chmidt integral
  operators and stochastic dynamical systems.
\newblock \emph{J. Nonlinear Sci.}, 32:\penalty0 79, 2022.
\newblock \doi{10.1007/s00332-022-09839-0}.

\bibitem[Santos~Guti{\'e}rez et~al.(2025)Santos~Guti{\'e}rez, Zagli, and
  Carigi]{SantosGutierezEtAl25}
M.~Santos~Guti{\'e}rez, N.~Zagli, and G.~Carigi.
\newblock Markov matrix perturbations to optimize dynamical and entropy
  functionals, 2025.
\newblock URL \url{https://arxiv.org/abs/2507.14040}.

\bibitem[Das and Giannakis(2019)]{DasGiannakis19}
S.~Das and D.~Giannakis.
\newblock Delay-coordinate maps and the spectra of {K}oopman operators.
\newblock \emph{J. Stat. Phys.}, 175\penalty0 (6):\penalty0 1107--1145, 2019.
\newblock \doi{10.1007/s10955-019-02272-w}.

\bibitem[Giannakis(2021)]{Giannakis21a}
D.~Giannakis.
\newblock Delay-coordinate maps, coherence, and approximate spectra of
  evolution operators.
\newblock \emph{Res. Math. Sci.}, 8:\penalty0 8, 2021.
\newblock \doi{10.1007/s40687-020-00239-y}.

\bibitem[Takens(1981)]{Takens81}
F.~Takens.
\newblock Detecting strange attractors in turbulence.
\newblock In \emph{Dynamical Systems and Turbulence}, volume 898 of
  \emph{Lecture Notes in Mathematics}, pages 366--381. Springer, Berlin, 1981.
\newblock \doi{10.1007/bfb0091924}.

\bibitem[Hameed et~al.(1999)Hameed, Goswami, Vinayachandran, and
  Yamagata]{HameedEtAl99}
S.~N. Hameed, B.~N. Goswami, P.~N. Vinayachandran, and T.~Yamagata.
\newblock A dipole mode in the tropical {I}ndian {O}cean.
\newblock \emph{Nature}, 401:\penalty0 360--363, 1999.
\newblock \doi{10.1038/43854}.

\bibitem[Webster et~al.(1999)Webster, Moore, Loschnigg, and
  Leban]{WebsterEtAl99}
P.~J. Webster, A.~Moore, J.~Loschnigg, and M.~Leban.
\newblock Coupled ocean dynamics in the {I}ndian {O}cean during 1997--98.
\newblock \emph{Nature}, 401:\penalty0 356--360, 1999.
\newblock \doi{10.1038/43848}.

\bibitem[Sauer et~al.(1991)Sauer, Yorke, and Casdagli]{SauerEtAl91}
T.~Sauer, J.~A. Yorke, and M.~Casdagli.
\newblock Embedology.
\newblock \emph{J. Stat. Phys.}, 65\penalty0 (3--4):\penalty0 579--616, 1991.
\newblock \doi{10.1007/bf01053745}.

\bibitem[Robinson(2005)]{Robinson05}
J.~C. Robinson.
\newblock A topological delay embedding theorem for infinite-dimensional
  dynamical systems.
\newblock \emph{Nonlinearity}, 18\penalty0 (5):\penalty0 2135--2143, 2005.
\newblock \doi{10.1088/0951-7715/18/5/013}.

\bibitem[Robinson(2008)]{Robinson08}
J.~C. Robinson.
\newblock A topological time-delay embedding theorem for infinite-dimensional
  cocycle dynamical systems.
\newblock \emph{Discrete Cont. Dyn. Syst. Ser. B}, 9\penalty0 (3\&4):\penalty0
  731--741, 2008.
\newblock \doi{10.3934/dcdsb.2008.9.731}.

\bibitem[Deyle and Sugihara(2011)]{DeyleSugihara11}
E.~R. Deyle and G.~Sugihara.
\newblock Generalized theorems for nonlinear state space reconstruction.
\newblock \emph{PLoS ONE}, 6\penalty0 (3):\penalty0 e18295, 2011.
\newblock \doi{10.1371/journal.pone.0018295}.

\bibitem[Nadaraya(1964)]{Nadaraya64}
E.~A. Nadaraya.
\newblock On estimating regression.
\newblock \emph{Theory Probab. Appl.}, 9\penalty0 (1):\penalty0 141--142, 1964.
\newblock \doi{10.1137/1109020}.

\bibitem[Watson(1964)]{Watson64}
G.~S. Watson.
\newblock Smooth regression analysis.
\newblock \emph{Sankhya Ser. A}, 26\penalty0 (4):\penalty0 359--372, 1964.

\bibitem[Stein(1993)]{Stein93}
E.~M. Stein.
\newblock \emph{Harmonic Analysis: Real-Variable Methods, Orthogonality, and
  Oscillatory Integrals}, volume~43 of \emph{Princeton Mathematical Series}.
\newblock Princeton University Press, Princeton, 1993.

\bibitem[Haj\l{}asz et~al.(2008)Haj\l{}asz, Koskela, and
  Tuominen]{hajlasz_sobolev_2008}
P.~Haj\l{}asz, P.~Koskela, and H.~Tuominen.
\newblock Sobolev embeddings, extensions and measure density condition.
\newblock \emph{Journal of Functional Analysis}, 254\penalty0 (5):\penalty0
  1217--1234, 2008.
\newblock ISSN 0022-1236.
\newblock \doi{10.1016/j.jfa.2007.11.020}.

\bibitem[Giannakis and Latifi-Jebelli(2026)]{GiannakisLatifiJebelli26}
D.~Giannakis and M.~J. Latifi-Jebelli.
\newblock Kernel smoothing operators on thick open domains.
\newblock \emph{Anal. Math.}, 2026.
\newblock In press.

\bibitem[Bj{\"o}rn et~al.(2001)Bj{\"o}rn, MacManus, and
  Shanmugalingam]{BjornEtAl01}
J.~Bj{\"o}rn, P.~MacManus, and N.~Shanmugalingam.
\newblock Fat sets and pointwise boundary estimates for $p$-harmonic functions
  in metric spaces.
\newblock \emph{J. Anal. Math.}, 85:\penalty0 339--369, 2001.
\newblock \doi{10.1007/bf02788087}.

\bibitem[Canto et~al.(2025)Canto, Ihnatsyeva, Lehrb{\"a}ck, and
  V{\"a}h{\"a}kangas]{CantoEtAl25}
J.~Canto, L.~Ihnatsyeva, J.~Lehrb{\"a}ck, and A.~V. V{\"a}h{\"a}kangas.
\newblock Capacities and density conditions in metric spaces.
\newblock \emph{Potential Anal.}, 62\penalty0 (2), 2025.
\newblock \doi{10.1007/s11118-024-10137-5}.

\bibitem[Khas'minskii(1963)]{Khasminskii63}
R.~Z. Khas'minskii.
\newblock The behavior of a self-oscillating system acted upon by slight noise.
\newblock \emph{Journal of Applied Mathematics and Mechanics}, 27\penalty0
  (4):\penalty0 1035--1044, 1963.
\newblock ISSN 0021-8928.
\newblock \doi{https://doi.org/10.1016/0021-8928(63)90184-9}.

\bibitem[Kifer(1974)]{Kifer74}
Y.~I. Kifer.
\newblock On small random perturbations of some smooth dynamical systems.
\newblock \emph{Math. USSR-Izv.}, 8:\penalty0 1083--1107, 1974.
\newblock \doi{10.1070/IM1974v008n05ABEH002139}.

\bibitem[Ulam(1964)]{Ulam64}
S.~M. Ulam.
\newblock \emph{Problems in Modern Mathematics}.
\newblock Dover Publications, Mineola, 1964.

\bibitem[Kemeny and Snell(1983)]{kemenysnell}
J.~G. Kemeny and J.~L. Snell.
\newblock \emph{Finite Markov Chains}.
\newblock Springer New York, 1983.

\bibitem[Froyland et~al.(2019)Froyland, Rock, and Sakellariou]{Froyland19SEBA}
G.~Froyland, C.~P. Rock, and K.~Sakellariou.
\newblock Sparse eigenbasis approximation: Multiple feature extraction across
  spatiotemporal scales with application to coherent set identification.
\newblock \emph{Communications in Nonlinear Science and Numerical Simulation},
  77:\penalty0 81--107, 2019.
\newblock \doi{https://doi.org/10.1016/j.cnsns.2019.04.012}.

\bibitem[Froyland et~al.(2024)Froyland, Giannakis, Luna, and
  Slawinska]{FroylandEtAl24}
G.~Froyland, D.~Giannakis, E.~Luna, and J.~Slawinska.
\newblock Revealing trends and persistent cycles of non-autonomous systems with
  operator-theoretic techniques: {A}pplications to past and present climate
  dynamics.
\newblock \emph{Nat. Commun.}, 15:\penalty0 4268, 2024.
\newblock \doi{10.1038/s41467-024-48033-6}.

\bibitem[Kato(1995)]{Kato95}
T.~Kato.
\newblock \emph{Perturbation Theory for Linear Operators}.
\newblock Classics in Mathematics. Springer-Verlag, Berlin, 2 edition, 1995.

\bibitem[Horn and Johnson(2013)]{HornJohnson13}
R.~A. Horn and C.~R. Johnson.
\newblock \emph{Matrix Analysis}.
\newblock Cambridge University Press, Cambridge; New York, 2nd edition, 2013.

\bibitem[Lorenz(1963)]{Lorenz63}
E.~N. Lorenz.
\newblock Deterministic nonperiodic flow.
\newblock \emph{J. Atmos. Sci.}, 20:\penalty0 130--141, 1963.
\newblock \doi{10.1175/1520-0469(1963)020<0130:dnf>2.0.co;2}.

\bibitem[Tucker(1999)]{Tucker99}
W.~Tucker.
\newblock The {L}orenz attractor exists.
\newblock \emph{C. R. Acad. Sci. Paris, Ser. I}, 328:\penalty0 1197--1202,
  1999.

\bibitem[Luzzatto et~al.(2005)Luzzatto, Melbourne, and Paccaut]{LuzzattoEtAl05}
S.~Luzzatto, I.~Melbourne, and F.~Paccaut.
\newblock The {L}orenz attractor is mixing.
\newblock \emph{Comm. Math. Phys.}, 260\penalty0 (2):\penalty0 393--401, 2005.

\bibitem[Froyland and Dellnitz(2003)]{FroylandDellnitz03}
G.~Froyland and M.~Dellnitz.
\newblock Detecting and locating near-optimal invariant sets and cycles.
\newblock \emph{SIAM J. Sci. Comput.}, 24\penalty0 (6):\penalty0 1839--1863,
  2003.
\newblock \doi{10.1137/s106482750238911x}.

\bibitem[Froyland and Padberg(2009)]{FroylandPadberg09}
G.~Froyland and K.~Padberg.
\newblock Almost-invariant sets and invariant manifolds -- {C}onnecting
  probabilistic and geometric descriptions of coherent structures in flows.
\newblock \emph{Phys. D}, 238:\penalty0 1507--1523, 2009.
\newblock \doi{10.1016/j.physd.2009.03.002}.

\bibitem[Beggs()]{RadialBasisFunctions_jl}
K.~Beggs.
\newblock Radialbasisfunctions.jl: Meshless {RBF} interpolation and
  differential operators for {Julia}.
\newblock URL \url{https://github.com/JuliaMeshless/RadialBasisFunctions.jl}.

\bibitem[Timmermann et~al.(2018)]{TimmermannEtAl18}
A.~Timmermann et~al.
\newblock El {N}i{\~n}o--{S}outhern {O}scillation complexity.
\newblock \emph{Nature}, 559:\penalty0 535--545, 2018.
\newblock \doi{10.1038/s41586-018-0252-6}.

\bibitem[Barnston et~al.(1997)Barnston, Chelliah, and
  Goldenberg]{BarnstonEtAl97}
A.~G. Barnston, M.~Chelliah, and S.~B. Goldenberg.
\newblock Documentation of a highly {ENSO}-related {SST} region in the
  equatorial {P}acific: {R}esearch note.
\newblock \emph{Atmos. Ocean}, 35\penalty0 (3):\penalty0 367--383, 1997.
\newblock \doi{10.1080/07055900.1997.9649597}.

\bibitem[Gent et~al.(2011)]{GentEtAl11}
P.~R. Gent et~al.
\newblock The {C}ommunity {C}limate {S}ystem {M}odel version 4.
\newblock \emph{J. Climate}, 24:\penalty0 4973--4991, 2011.
\newblock \doi{10.1175/2011jcli4083.1}.

\bibitem[Saji et~al.(1999)Saji, Goswami, Vinayachandran, and
  Yamagata]{SajiEtAl99}
H.~N. Saji, B.~N. Goswami, P.~N. Vinayachandran, and T.~Yamagata.
\newblock A dipole mode in the tropical {I}ndian {O}cean.
\newblock \emph{Nature}, 401:\penalty0 360--363, 1999.
\newblock \doi{10.1038/43854}.

\bibitem[Slawinska and Giannakis(2017)]{SlawinskaGiannakis17}
J.~Slawinska and D.~Giannakis.
\newblock Indo-{P}acific variability on seasonal to multidecadal time scales.
  {P}art {I}: {I}ntrinsic {SST} modes in models and observations.
\newblock \emph{J. Climate}, 30\penalty0 (14):\penalty0 5265--5294, 2017.
\newblock \doi{10.1175/jcli-d-16-0176.1}.

\bibitem[Hameed et~al.(2018)Hameed, Jin, and Thilakan]{HameedEtAl18}
S.~N. Hameed, D.~Jin, and V.~Thilakan.
\newblock A model for super {E}l {N}i\~{n}os.
\newblock \emph{Nature Communications}, 2018.
\newblock \doi{10.1038/s41467-018-04803-7}.

\end{thebibliography}

\end{document}